\documentclass{siamltex}

\usepackage[breaklinks,colorlinks=true,linkcolor=blue,citecolor=red,backref=page]{hyperref}

\usepackage{amsfonts} 
\usepackage{amssymb,amsmath} 
\usepackage{graphicx,color}

\DeclareMathAlphabet{\itbf}{OML}{cmm}{b}{it}

\def\bb{{{\itbf b}}}
\def\bq{{{\itbf q}}}

\def\br{{{\itbf r}}}

\def\by{{{\itbf y}}}
\def\bx{{{\itbf x}}}

\def\bh{{{\itbf h}}}

\def\bk{{{\itbf k}}}

\def\bzeta{{\boldsymbol{\zeta}}}
\def\bxi{{\boldsymbol{\xi}}}
\def\balpha{{\boldsymbol{\alpha}}}
\def\bbeta{{\boldsymbol{\beta}}}

\def\eps{{\varepsilon}}

\newcommand{\RR}{\mathbb{R}}

\newcommand{\EE}{\mathbb{E}}

\newtheorem{thm}{Theorem}[section]

\newtheorem{prop}[thm]{Proposition}

\begin{document}  
 
 \title{Focusing Waves Through a Randomly Scattering Medium in the White-Noise Paraxial Regime} 

\author{Josselin Garnier\footnotemark[1] 
\and Knut S\O lna\footnotemark[2]}  

\maketitle

\footnotetext[1]{Laboratoire de Probabilit\'es et Mod\`eles Al\'eatoires
\& Laboratoire Jacques-Louis Lions,
Universit{\'e} Paris Diderot,
75205 Paris Cedex 13,
France
{\tt garnier@math.univ-paris-diderot.fr}}
\footnotetext[2]{Department of Mathematics, 
University of California, Irvine CA 92697
{\tt ksolna@math.uci.edu}}

\renewcommand{\thefootnote}{\arabic{footnote}}

\begin{abstract}
When waves propagate through a complex or heterogeneous medium the  wave field is 
corrupted by the heterogeneities. Such corruption limits the performance of imaging or
communication schemes.  One may then ask the question: is there an optimal way of encoding
a signal so as to counteract the corruption by the medium? In the ideal situation
the answer is given by  time reversal: for a given target or focusing point, in a first step
let the target emit a signal and then record the signal transmitted to the source antenna, 
time reverse this and use it as the  source trace at the source antenna in a second step.
This source will give a sharply focused wave at the target location if the source aperture is large enough.
Here we address this scheme in the more practical situation with a limited aperture,
time-harmonic signal, and finite-sized elements in the source array.  Central questions 
are  then the focusing resolution and signal-to-noise ratio at the target,
their  dependence on the physical parameters, and the capacity to focus selectively in the
neighborhood of the target point and therefore to transmit images. Sharp results are presented for  these questions. 

\end{abstract}



\pagestyle{myheadings}
\thispagestyle{plain}
\markboth{J. Garnier and K. S\o lna}
{Focusing Waves Through a Randomly Scattering Medium}

\section{Introduction}

\begin{figure}
\begin{center}
\begin{picture}(330,160)
\put(0,117){{\bf (a)}}
\put(10,80){\includegraphics[width=7.7cm]{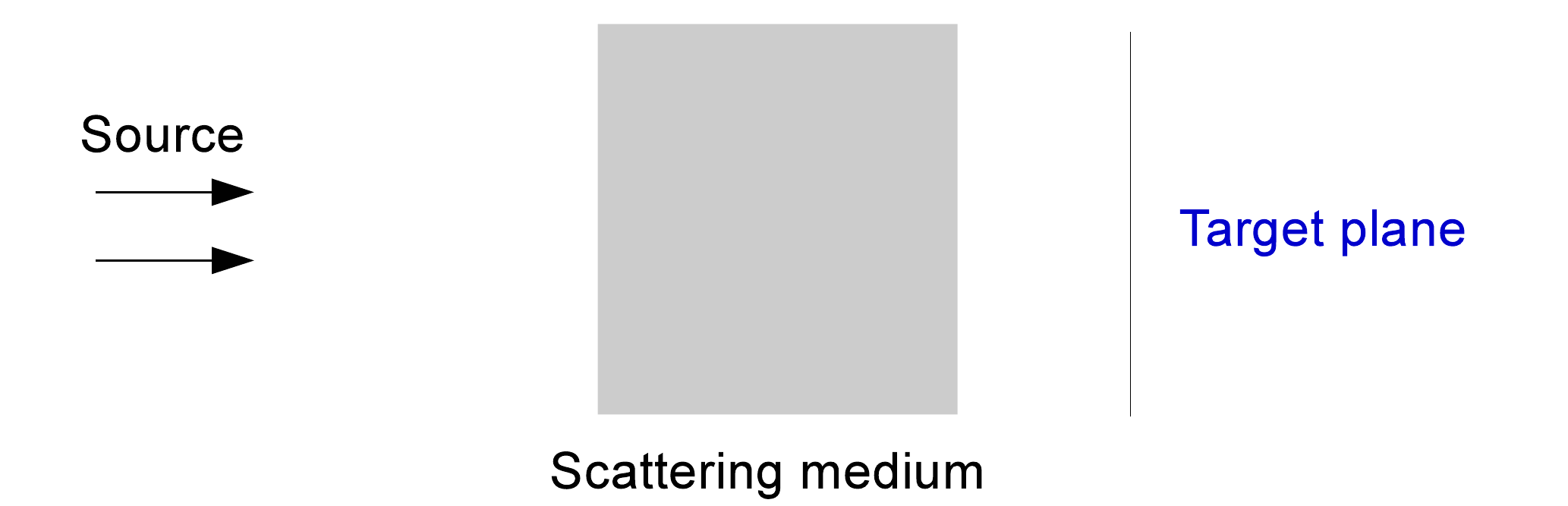}}
\put(225,84){\includegraphics[width=3.3cm]{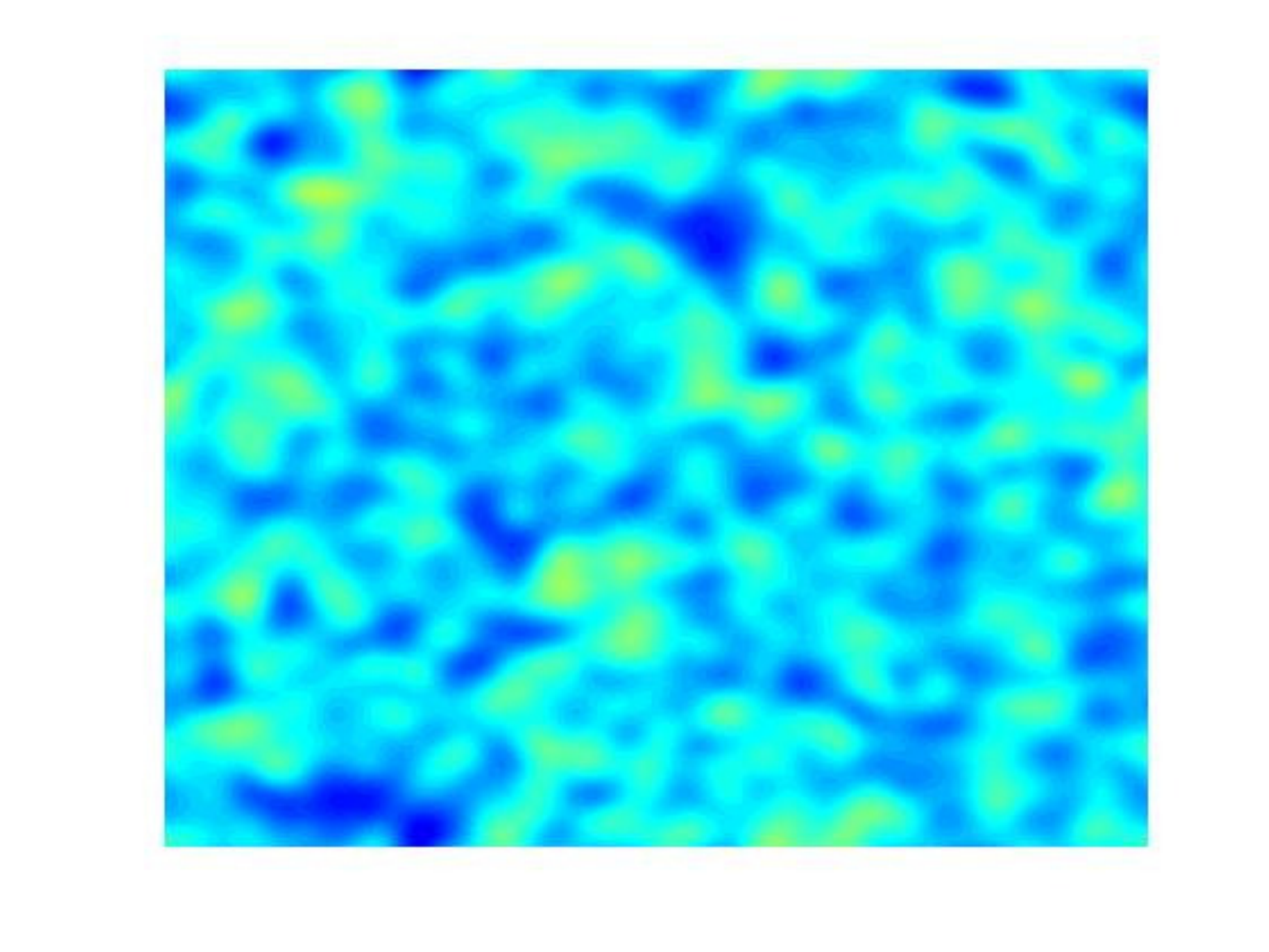}}
\put(0,37){{\bf (b)}} 
\put(10,0){\includegraphics[width=7.7cm]{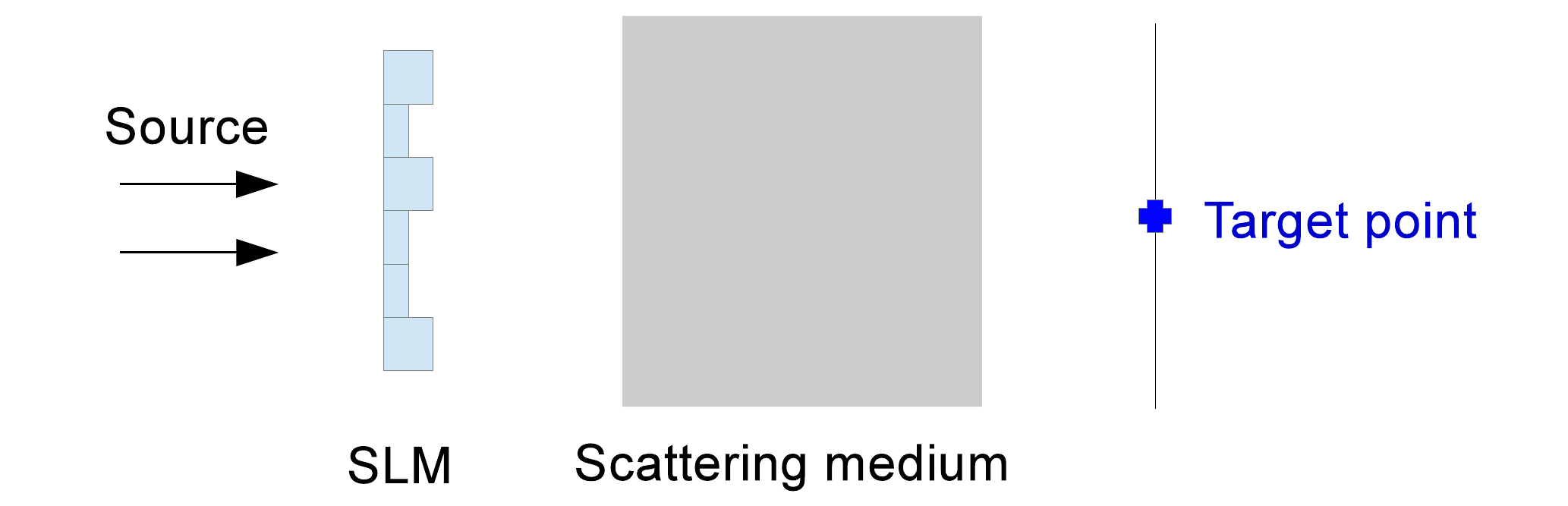}}
\put(225,4){\includegraphics[width=3.3cm]{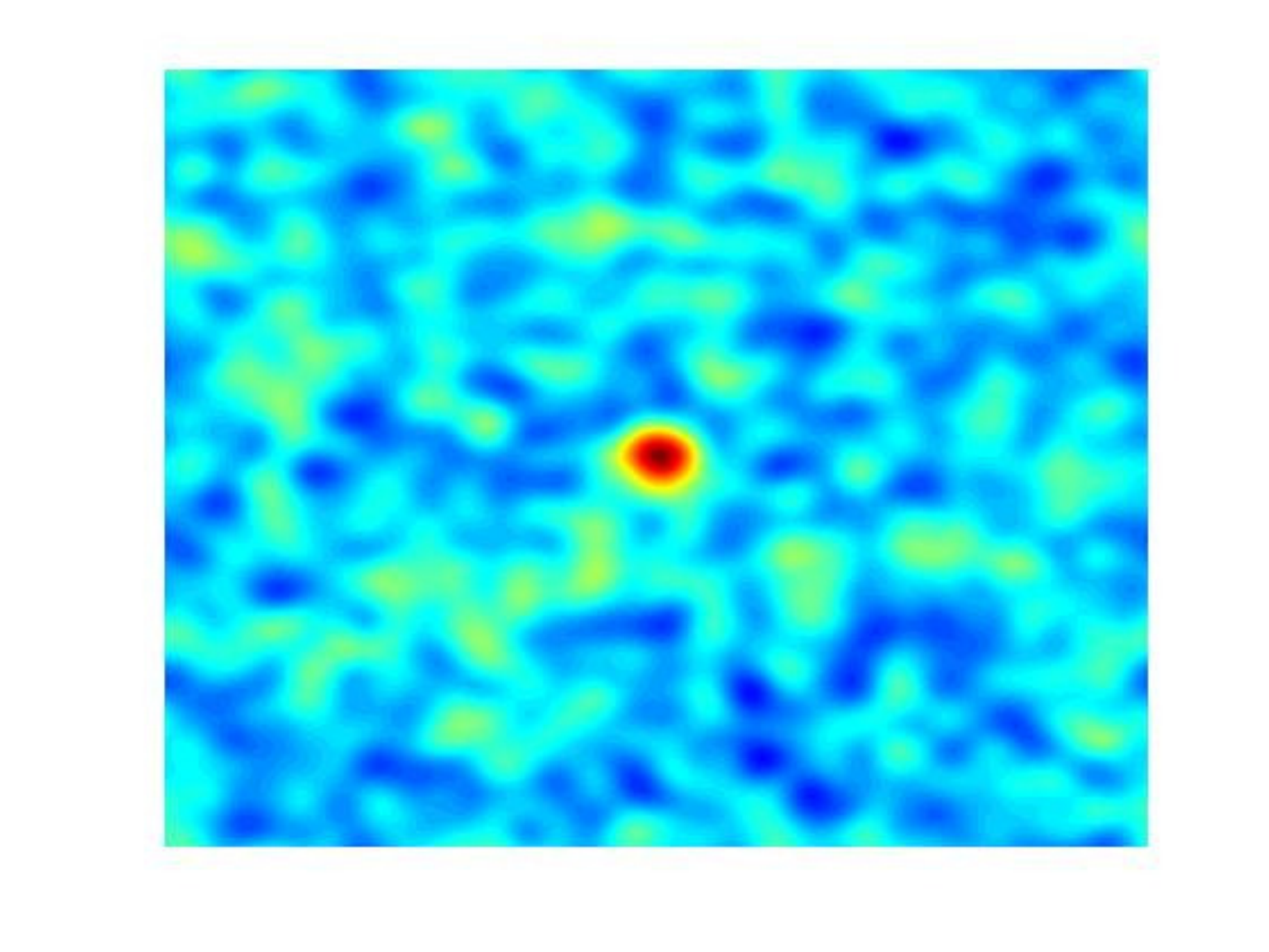}}
\end{picture}
\end{center}
\caption{Focusing wave through a scattering medium. Without any control 
one gets a speckle pattern in the target plane (a).
With a spatial light modulator (SLM) one can focus on a target point by imposing appropriate phase shifts (b).
\label{fig:0a} 
}
\end{figure}

Wavefront-shaping-based schemes for
focusing \cite{popoff14,vellekoop10,vellekoop07,vellekoop08}  and
 imaging \cite{popoff10,mosk12,katz12}
have proved very useful for focusing and imaging through scattering media.
The primary goal of the experiments reported in these papers
is to focus monochromatic light through a layer of strongly scattering material.
This is a challenging problem as it is known that multiple scattering of waves  by the medium inhomogeneities
scrambles the transmitted light into random interference patterns called speckle patterns \cite{goodman}.
However, if a spatial light modulator (SLM) is applied before the scattering medium, then it is possible
to focus light as first demonstrated in \cite{vellekoop07}. 
Indeed, the elements of the SLM can impose phase shifts prescribed by the user, and it is possible to choose 
(by an optimization scheme) the phase shifts so as to maximize the 
intensity transmitted at one point in the target plane behind the scattering medium (see Figure \ref{fig:0a}).

\begin{figure}
\begin{center}
\begin{picture}(330,160)
\put(0,117){{\bf (a)}}
\put(10,80){\includegraphics[width=7.7cm]{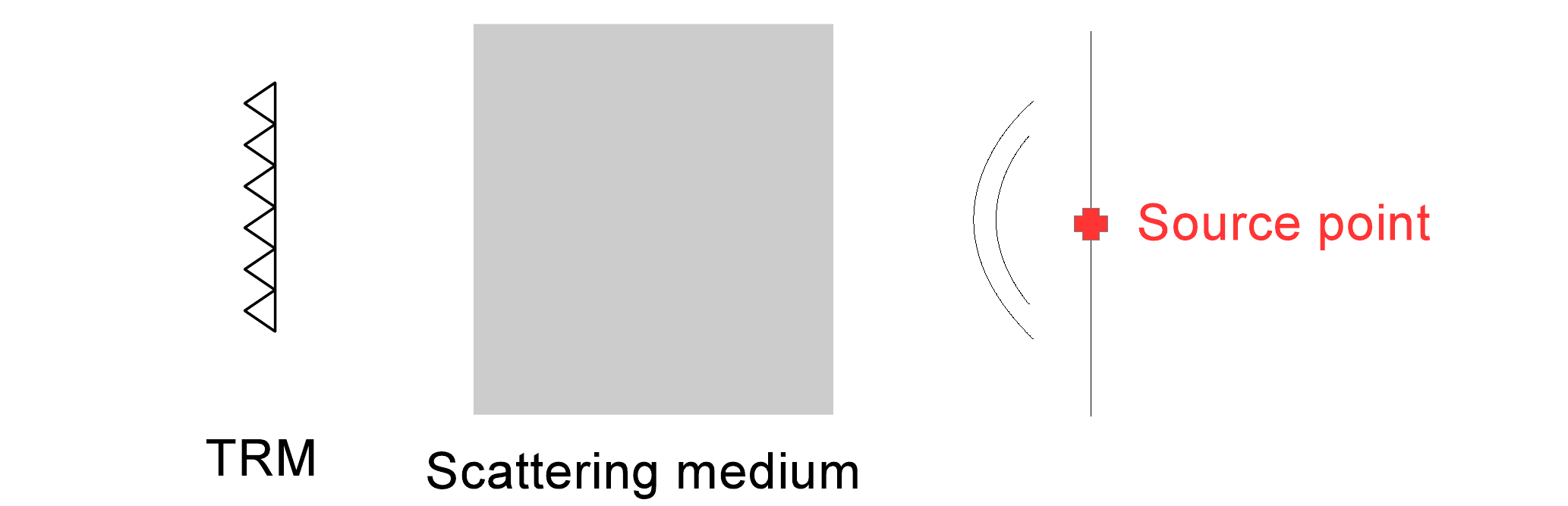}}
\put(0,37){{\bf (b)}} 
\put(10,0){\includegraphics[width=7.7cm]{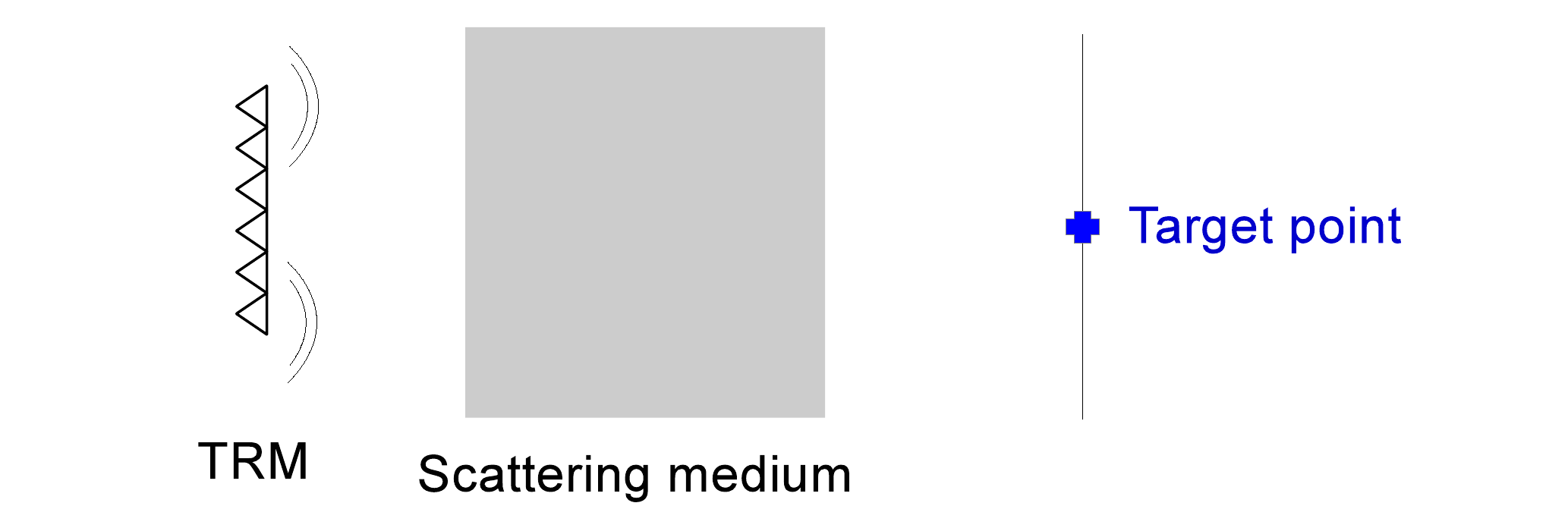}}
\put(205,-4){ \includegraphics[width=4.0cm]{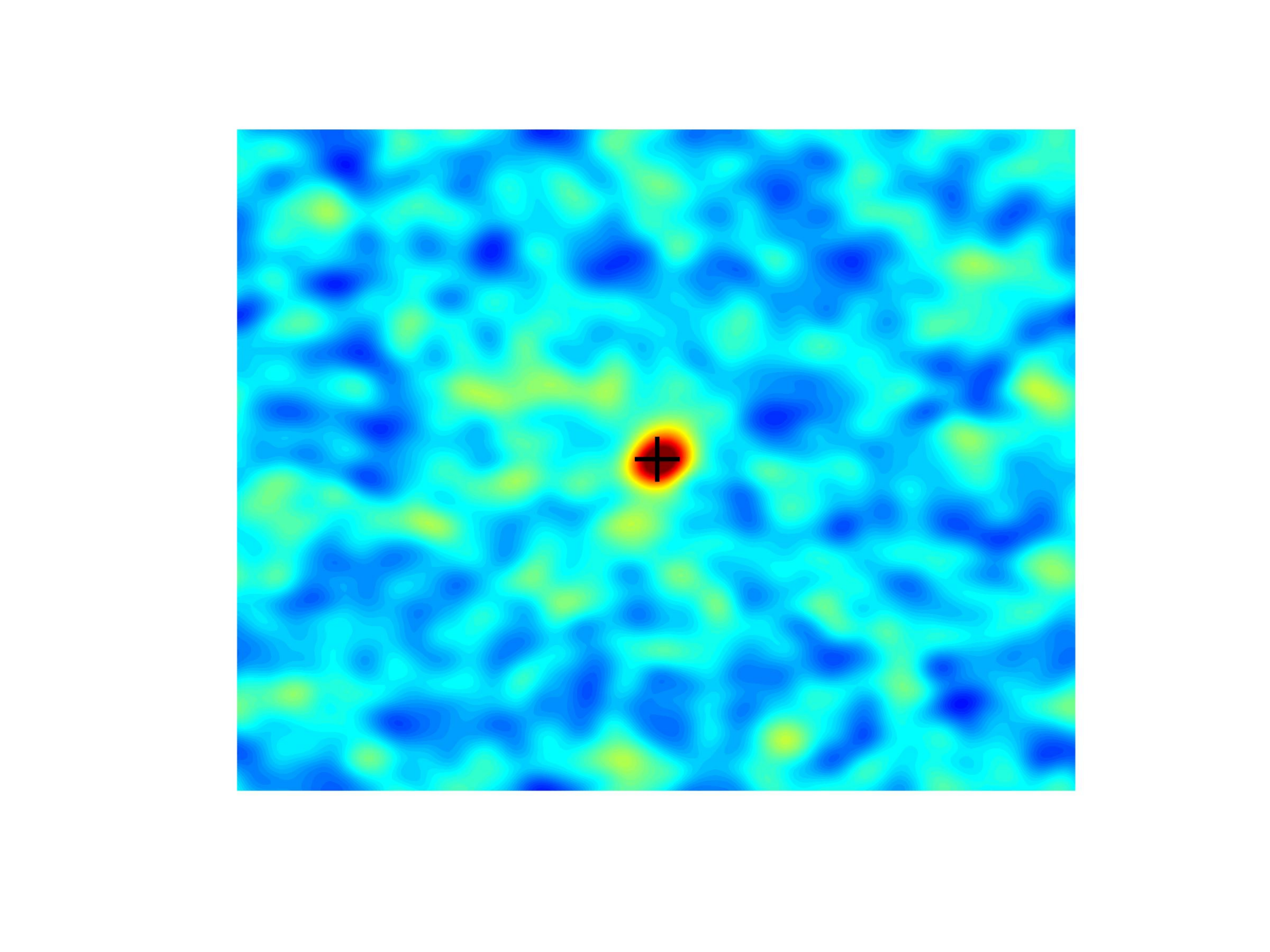}}
\end{picture}
\end{center}
\caption{Time-reversal experiment through a scattering medium. In the first step of the experiment (a)
a time-harmonic point source emits a wave that propagates through the scattering medium and is recorded by the time-reversal mirror (TRM) used as an array of receivers. 
In the second step of the experiment (b) the TRM is used as an array of sources, it
emits the complex-conjugated recorded field, and the wave refocuses at the original source location (the cross in the right image
stands for the original source location).
\label{fig:0b} 
}
\end{figure}

It turns out that the phase shifts obtained by the wavefront-shaping optimization procedure 
are the opposite phases of the field emitted by a point source at the target point and recorded in the plane of the SLM \cite{mosk12}.
In other words, the wavefront-shaping optimization procedure is equivalent to 
phase conjugation or time reversal.
The focal spot obtained at the target point by the wavefront-shaping-based scheme
is the focal spot of the time-reversed refocused wave obtained at the end of a time-reversal experiment 
in which waves emitted by a source at the target point propagate through the scattering medium,
are recorded by an array of sensors in the plane of the SLM, are time-reversed and re-emitted 
through the scattering medium towards the plane of the target  (see Figure \ref{fig:0b}).
This time-reversal interpretation and the known refocusing and stability properties 
of time reversal for waves in random media \cite{lerosey07,book1} explain the focusing properties
of the wavefront-shaping-based scheme \cite{mosk12}.

However, two questions can be raised in view of the experimental results:

-  the elements of the SLM are much larger than the operating wavelength 
(for instance, $20\mu$m for the SLM used in \cite{vellekoop07,vellekoop08})
and may even be larger than the correlation radius of a field emitted from the target point
and transmitted through the scattering medium.
This means that the phase shifts imposed by the elements of the SLM are not the phases of the conjugated
(or time-reversed) field,
but result from local averages of this field.
Nevertheless refocusing seems very efficient. We will see that indeed time reversal with a 
smoothing kernel whose radius is larger than the correlation radius allows for wave refocusing, 
that the radius of the focal spot at the target point does not  depend on the smoothing, 
but that the statistical stability of the focal spot (i.e. the signal-to-noise ratio)  depends on it.

- the SLM is used to correct the wavefront distortions of the diffused light. 
After performing the correction for a target point, it turns out that the correction can also be used to 
make the wave focus on a prescribed neighboring point, which allows for imaging of objects. This is explained in the physical literature by 
the memory effect for speckle correlations \cite{feng88,freund88,katz12}.
We will see that indeed it is possible to focus selectively at any point in a neighborhood 
of the target point, which allows to image an object, and we will quantify the 
extent of this neighborhood and the resolution and the stability of the image obtained by this method.

We address in this paper these questions in the paraxial white-noise regime,
as described by the  It\^o-Schr\"odinger model for the Green's function.
This model is  a simplification of the Helmholtz equation with random index of refraction, 
it gives the correct statistical structure of the wave field when the propagation distance
is larger than the correlation length of the medium which is itself larger than the wavelength and
when the typical amplitude of the medium fluctuations is small.
The It\^o-Schr\"odinger model  can be derived  rigorously from  the Helmholtz equation 
by a separation of scales technique in the high-frequency regime  \cite{garniers0,garniers1,garniers2}.
It is physically relevant and it models  many situations, for instance laser beam propagation \cite{andrews,strohbehn},
time reversal in random media  \cite{blomgren,PRS04}, or
underwater acoustics \cite{tappert}.
The It\^o-Schr\"odinger model allows for the use of  It\^o's stochastic calculus,  
which in turn enables the  closure of the hierarchy of moment equations \cite{fps,ishimaru}.
Until recently, the equation for the fourth-order moments of the Green's function
could not be solved \cite[Sec.~20.18]{ishimaru}.
However in a recent paper \cite{garniers4} (with a preliminary version in \cite{garniers3})
the behavior of the fourth-order moments of the random 
 paraxial Green's function could be unraveled and our paper is based 
 on this result that allows to carry out a variance analysis of the time-reversed field
 (which indeed involves fourth-order moments of the Green's function).
 
The paper is organized as follows. 
We describe the transmission problem in Section \ref{sec:1}.
In Section \ref{sec:regime} we extend the result obtained in \cite{garniers4} 
in order to get the fourth-order moment that is of interest for our study.
In Section \ref{sec:refocus} we prove Propositions \ref{prop:momtr1} and \ref{prop:momtr2}
that describe the mean and covariance of the refocused field around the target point.
In Section \ref{sec:prescribed} we show how to focus on a prescribed point in the neighborhood of the 
target point and use this idea in Section \ref{sec:image} to show that it possible to transmit an image through a strongly
scattering medium.


\section{Time-Reversal Experiment}
\label{sec:1}%
In this paper we denote the three-dimensional spatial variable by $(\bx,z)$, with $\bx \in \RR^2$ the transverse
variable and $z\in \RR$ the longitudinal variable.
A time-reversal mirror (TRM) is located in the plane $z=0$. Its radius is $R_{\rm m}$ and the radius of its elements is $\rho_0$.

In the first step of the time-reversal experiment, 
a point  source localized at $(\by,L)$ emits a time-harmonic signal at frequency $\omega$.
The TRM is used as an array of receivers 
and records the wave emitted by the point source. The size $\rho_0$ of the elements of the TRM is taken into account 
in the form of a Gaussian smoothing kernel with radius $\rho_0$. 
By denoting $\hat{\cal G}(L , \bx_m,\by)$ the Green's function 
from $(\bx_m,0)$ to $(\by,L)$ (which is equal to the Green's function from $(\by,L)$
to $(\bx_m,0)$ by reciprocity), the recorded field at $(\bx_m,0)$ is  therefore
obtained  via the smoothing over  the Gaussian  mirror element as:  
\begin{equation}
\label{def:urec}
\hat{u}_{\rm rec}(\bx_m;\by) =
\frac{1}{2 \pi \rho_0^2} \int \hat{\cal G}(L , \bx_m+\bx',\by) \exp \Big( - \frac{|\bx'|^2}{2 \rho_0^2} \Big) d\bx'  .
\end{equation}

In the second step of the  time-reversal experiment, 
the TRM is used as an array of sources. It emits
the time-reversed (or complex-conjugated) recorded field $\overline{\hat{u}_{\rm rec}}$.
The time-reversed field observed in the plane $z=L$ at the point $(\bx,L)$ has the form
\begin{equation}
\label{eq:trfield1}
\hat{u}_{\rm tr}( \bx ;\by ) =
 \int \hat{u}_{\rm em} (\bx,\bx_m) \exp \Big( - \frac{|\bx_m|^2}{R_{\rm m}^2} \Big) 
\overline{ \hat{u}_{\rm rec}(\bx_m; \by) }  
  d\bx_m .
\end{equation}
Here we have assumed that the TRM has a radius $R_{\rm m}$ and
can be modeled by a Gaussian spatial cut-off function.  Moreover, 
we again take into account the size $\rho_0$ of the elements of the TRM by considering that
from any point $(\bx_m,0)$ the  TRM  can emit from a patch with radius $\rho_0$ and 
with a Gaussian form,
which generates the following field at point $(\bx,L)$:
\begin{equation}
\hat{u}_{\rm em} (\bx; \bx_m) = \frac{1}{2 \pi \rho_0^2} 
\int \hat{\cal G}(L,\bx_m+\bx',\bx) \exp \Big( - \frac{|\bx'|^2}{2 \rho_0^2} \Big) d\bx'  .
\end{equation}

The  time-reversed field observed in the plane $z=L$ can therefore be expressed as
\begin{equation}
\label{eq:base}
\hat{u}_{\rm tr}( \bx ;\by ) = 4k_0^2 C_0
\iint 
\exp \Big( - \frac{|\bx'|^2}{r_0^2} - \frac{|\by'|^2}{4 \rho_0^2} \Big) 
\hat{\cal G}\big(L,\bx'+\frac{\by'}{2},\bx\big) 
\overline{\hat{\cal G}\big(L,\bx'-\frac{\by'}{2},\by\big) }
  d\bx' d\by'   ,
\end{equation}
with   
\begin{equation}
\label{def:r0}
C_0 = \frac{r_0^2 -\rho_0^2}{16 \pi  k_0^2 \rho_0^2 r_0^2}
,\quad \quad 
r_0^2 = R_{\rm m}^2 + \rho_0^2  .
\end{equation}  
Of course, when the size of the elements $\rho_0$ goes to zero, we recover the standard expression
for the refocused time-reversed field with a time-reversal mirror with a Gaussian aperture with radius $R_{\rm m}$:
$$
\hat{u}_{\rm tr}( \bx ;\by ) \mid_{\rho \to 0}=  
\int 
\exp \Big( - \frac{|\bx'|^2}{R_{\rm m}^2}  \Big) 
\hat{\cal G}\big(L,\bx' ,\bx\big) 
\overline{\hat{\cal G}\big(L,\bx' ,\by\big) }
  d\bx'     .
$$
From now on we will take $C_0=1$ as this multiplicative factor does not play any role in what follows.

{\bf Remark.}
We have modeled the global shape of the TRM and the local shape of the elements of the TRM by soft Gaussian cut-off functions,
instead of hard cut-off functions such as ${\bf 1}_{[0,R_{\rm m}]}(|\bx_m|)$ or ${\bf 1}_{[0,\rho_0]}(|\bx'|)$,
 because this allows to get simple and explicit expressions in the following.

\subsection{The Green's Function in the White-noise Paraxial Regime}
In the white-noise paraxial regime the Green's function $\hat{\cal G}$ is of the form \cite{garniers4}
$$
\hat{\cal G} (L,\bx, \by\big)  = \frac{i}{2k_0} e^{i k_0 L}  \hat{G}(L,\bx,\by) ,
$$
where $k_0$ is the homogeneous wavenumber and 
 the function $\hat{G}$ is the solution of the It\^o-Schr\"odinger equation
\begin{equation}
\label{eq:model}
 d \hat{G}(z,\bx,\by)   =     
          \frac{ i }{2k_0} \Delta_{\bx}   \hat{G}(z,\bx,\by) dz
   +   \frac{ik_0}{2}   \hat{G} (z,\bx,\by ) \circ  d{B}(z,\bx) 
  , 
\end{equation}
with the initial condition in the plane $z=0$:
$ \hat{G} (z= 0,\bx,\by )  = \delta(\bx-\by)$.
Here the symbol $\circ$ stands for the Stratonovich stochastic integral,
  $B(z,\bx)$ is a real-valued  Brownian field over $[0,\infty) \times \RR^2$ with  covariance
 \begin{equation}
 \label{defB}
\EE[   {B}(z,\bx)  {B}(z',\bx') ] =  
 {\min\{z, z'\}}   {C}(\bx - \bx')   ,
\end{equation}
and $C$ is determined by the two-point statistics of the fluctuations of the random medium
(in particular the width of $C$ is the correlation length of the medium fluctuations).
In It\^o's form Eq.~(\ref{eq:model}) reads
\begin{equation}
\label{eq:modelito}
 d \hat{G}(z,\bx,\by)   =     
          \frac{ i }{2k_0} \Delta_{\bx}   \hat{G}(z,\bx,\by) dz
   +   \frac{ik_0}{2}   \hat{G} (z,\bx,\by )   d{B}(z,\bx) - \frac{k_0^2}{8} C({\bf 0}) \hat{G} (z,\bx,\by ) dz .
\end{equation}
This equation was analyzed for the first time in \cite{dawson84} and it was derived from first principles
by a multiscale analysis of the wave equation in a random medium in \cite{garniers1}.
In this context, the time-reversed field is
\begin{equation}
\hat{u}_{\rm tr}( \bx ;\by ) =  
\iint 
\exp \Big( - \frac{|\bx'|^2}{r_0^2} - \frac{|\by'|^2}{4 \rho_0^2} \Big) 
\hat{G}\big(L,\bx'+\frac{\by'}{2},\bx\big) 
\overline{\hat{G}\big(L,\bx'-\frac{\by'}{2},\by\big) }
  d\bx' d\by'   .
\end{equation}

\subsection{The Mean Refocused Wave}
The mean time-reversed field observed at $(\bx,L)$ when the original source is at $(\by,L)$ is
\begin{equation}
M_1(L,\bx,\by) = \EE \big[ \hat{u}_{\rm tr}( \bx ;\by )\big]   ,
\end{equation}
and it can be expressed as
$$
M_1(L,\bx,\by)=
\iint 
\exp \Big( - \frac{|\bx'|^2}{r_0^2} - \frac{|\by'|^2}{4 \rho_0^2} \Big) 
\EE \Big[ \hat{G}\big(L,\bx'+\frac{\by'}{2},\bx\big) 
\overline{\hat{G}\big(L,\bx'-\frac{\by'}{2},\by\big) } \Big]
  d\bx' d\by'   .
$$
$M_1$ satisfies the system:
\begin{equation}
\label{eq:mom11}
\frac{\partial M_{1} }{\partial z} = \frac{i}{2k_0}  \big(  \Delta_{\bx}
-  \Delta_{\by} \big) M_{1}  + \frac{k_0^2}{4} \big( C(\bx-\by)-C({\bf 0})  \big) M_{1}    ,
\end{equation}
starting from 
$$
 M_{1} (z=0,\bx,\by) =
 \exp \Big( - \frac{|\bx+\by|^2}{4 r_0^2} - \frac{|\bx-\by|^2}{ 4 \rho_0^2} \Big)  .
 $$
 After parameterizing the two points $\bx$ and $\by$ as
 $$
 \br= \frac{\bx+\by}{2} , \quad \quad \bq = \bx-\by, 
 $$
 this equation can be solved (after Fourier transforming in $\br$):
 \begin{eqnarray}
 \nonumber
 M_1 \big( L,\br+ \frac{\bq}{2},\br- \frac{\bq}{2}\big) =
 \frac{r_0^2}{4 \pi} 
 \int \exp \Big( i \bxi \cdot \br
 - \frac{r_0^2 |\bxi|^2}{4} -\frac{ |\bq - \bxi \frac{L}{k_0} |^2}{4 \rho_0^2} \\
 +\frac{k_0^2}{4} \int_0^L C ( \bq -\bxi \frac{z}{k_0} ) -C({\bf 0}) dz \Big) d\bxi .
 \label{eq:mean1}
 \end{eqnarray}

\subsection{The Fluctuations of the Refocused Wave}
We consider the general second-order moment for the time-reversed wave
\begin{eqnarray}
\nonumber
M_2(L,\bx_1,\bx_2,\by_1,\by_2) &=& \EE \big[ \hat{u}_{\rm tr}( \bx_1 ;\by_1 )
\overline{\hat{u}_{\rm tr}( \by_2 ;\bx_2 )} \big] \\
&=& 
 \EE \big[ \hat{u}_{\rm tr}( \bx_1 ;\by_1 )
{\hat{u}_{\rm tr}( \bx_2 ;\by_2 )} \big]   ,
\end{eqnarray}
which depends on the fourth-order moment of the paraxial Green's function:
\begin{eqnarray*}
\nonumber
 M_2(L,\bx_1,\bx_2,\by_1,\by_2) &=&
\iint 
\exp \Big( - \frac{|\bx_1'|^2+|\bx_2'|^2}{r_0^2} - \frac{|\by_1'|^2+|\by_2'|^2}{4 \rho_0^2} \Big) \\
&& 
\times \EE \Big[ \hat{G}\big(L,\bx_1'+\frac{\by_1'}{2},\bx_1\big) \hat{G}\big(L,\bx_2'+\frac{\by_2'}{2},\bx_2\big) \\
&& \quad\quad \times
\overline{\hat{G}\big(L,\bx_1'-\frac{\by_1'}{2},\by_1\big) } \overline{\hat{G}\big(L,\bx_2'-\frac{\by_2'}{2},\by_2\big) } \Big]
  d\bx_1' d\by_1' 
  d\bx_2' d\by_2' .
\end{eqnarray*}
It satisfies
\begin{eqnarray}
\label{def:generalmoment}
\frac{\partial M_2}{\partial z} = \frac{i}{2k_0}  \Big(  \Delta_{\bx_1}+  \Delta_{\bx_2}
-   \Delta_{\by_1}- \Delta_{\by_2} \Big) M_{2} + \frac{k_0^2}{4} U_{2} \big( \bx_1,\bx_2,\by_1,\by_2 \big)
 M_{2}  ,
\end{eqnarray}
with the generalized potential
\begin{eqnarray}
\nonumber
 U_{2}\big(  \bx_1,\bx_2,\by_1,\by_2 \big)  &=&
 {C}(\bx_1-\by_1) + {C}(\bx_1-\by_2) + {C}(\bx_2-\by_1) + {C}(\bx_2-\by_2) \\
&&-  {C}( \bx_1-\bx_2)
- {C}( \by_1-\by_2) -
2{C}({\bf 0}) \, ,
\end{eqnarray}
and it starts from
\begin{eqnarray*}
 M_{2} (z=0,  \bx_1,\bx_2,\by_1,\by_2) &=&
 \exp \Big( - \frac{|\bx_1+\by_1|^2+|\bx_2+\by_2|^2}{4 r_0^2} \Big)\\
 && \times \exp \Big(- \frac{|\bx_1-\by_1|^2+|\bx_2-\by_2|^2}{ 4 \rho_0^2} \Big) .
\end{eqnarray*}

We parameterize  the four points 
$\bx_1,\bx_2,\by_1,\by_2$ in (\ref{def:generalmoment}) in the special way:
\begin{eqnarray*}
\bx_1 = \frac{\br_1+\br_2+\bq_1+\bq_2}{2}, \quad \quad 
\by_1 = \frac{\br_1+\br_2-\bq_1-\bq_2}{2}, \\
\bx_2 = \frac{\br_1-\br_2+\bq_1-\bq_2}{2}, \quad \quad 
\by_2 = \frac{\br_1-\br_2-\bq_1+\bq_2}{2}.
\end{eqnarray*}
In particular $\br_1/2$ is the barycenter of the four points $\bx_1,\bx_2,\by_1,\by_2$:
\begin{eqnarray*}
\br_1 = \frac{\bx_1+\bx_2+\by_1+\by_2}{2} , \quad \quad 
\bq_1 = \frac{\bx_1+\bx_2-\by_1-\by_2}{2}, \\
\br_2 = \frac{\bx_1-\bx_2+\by_1-\by_2}{2}, \quad \quad 
\bq_2 = \frac{\bx_1-\bx_2-\by_1+\by_2}{2}.
\end{eqnarray*}
In the variables $(\bq_1,\bq_2,\br_1,\br_2)$ the function $M_{2}$ satisfies the equation:
\begin{equation}
\label{eq:M20}
\frac{\partial M_{2}}{\partial z} = \frac{i}{k_0} \big( \nabla_{\br_1}\cdot \nabla_{\bq_1}
+
 \nabla_{\br_2}\cdot \nabla_{\bq_2}
\big)  M_{2} + \frac{k_0^2}{4} U_{2}(\bq_1,\bq_2,\br_1,\br_2) M_{2}   ,
\end{equation}
with the generalized potential
\begin{eqnarray}
\nonumber
U_{2}(\bq_1,\bq_2,\br_1,\br_2) &=& 
{C}(\bq_2+\bq_1) 
+
{C}(\bq_2-\bq_1) 
+
{C}(\br_2+\bq_1) 
+
{C}(\br_2-\bq_1) \\
&&- 
{C}( \bq_2+\br_2) - {C}( \bq_2-\br_2) - 2 {C}({\bf 0})  .
\end{eqnarray}
Note in particular that the generalized potential does not depend on the barycenter $\br_1$ as the medium
is statistically homogeneous.
The initial condition for Eq.~(\ref{eq:M20}) is
$$
M_{2}(z=0,\bq_1,\bq_2,\br_1,\br_2) = \exp \Big( - \frac{  |\br_1|^2+ |\br_2|^2}{2r_0^2}
- \frac{ |\bq_1|^2+ |\bq_2|^2}{2\rho_0^2} \Big).
$$
The second moment of the refocused wave field can be expressed as:
\begin{equation}
\label{eq:secondmom1}
 \EE \big[ \big| \hat{u}_{\rm tr} (\bx , \by)\big|^2  \big] = M_2( L, \br_1=\bx+\by, \br_2= {\bf 0} , \bq_1={\bf 0}, \bq_2 = \bx-\by ),
\end{equation}
and, more generally, 
\begin{equation}
\label{eq:secondmom1b}
 \EE \Big[  \hat{u}_{\rm tr} \big(\bx +\frac{\bh}{2} , \by\big)  
 \overline{ \hat{u}_{\rm tr} \big(\bx-\frac{\bh}{2} , \by\big)}  \Big] = M_2
 \big( L, \br_1=\bx+\by, \br_2= \frac{\bh}{2} , \bq_1= \frac{\bh}{2} , \bq_2 =  \bx-\by \big) .
\end{equation}

The Fourier transform (in $\bq_1$, $\bq_2$, $\br_1$, and $\br_2$) of the fourth-order moment
of the paraxial Green's function 
is defined by:
\begin{eqnarray}
\nonumber
\hat{M}_{2}(z,\bxi_1,\bxi_2,\bzeta_1,\bzeta_2) 
&=& 
\iint M_{2}(z,\bq_1,\bq_2,\br_1,\br_2)  \\
&& \hspace*{-0.4in}
\times
\exp  \big(- i\bq_1 \cdot \bxi_1- i\br_1\cdot \bzeta_1- i\bq_2 \cdot \bxi_2- i\br_2\cdot \bzeta_2\big) d\br_1d\br_2 d\bq_1d\bq_2
. \hspace*{0.3in} \end{eqnarray}
It satisfies
\begin{eqnarray}
\nonumber
&& 
\frac{\partial \hat{M}_{2}}{\partial z} + \frac{i}{k_0} \big( \bxi_1\cdot \bzeta_1+   \bxi_2\cdot \bzeta_2\big) \hat{M}_{2}
=
\frac{k_0^2}{4 (2\pi)^2} 
\int \hat{C}(\bk) \bigg[  
 \hat{M}_{2} (  \bxi_1-\bk, \bxi_2-\bk, \bzeta_1, \bzeta_2)  \\
\nonumber
&& \hspace*{0.7in} + 
 \hat{M}_{2} (  \bxi_1-\bk,\bxi_2,  \bzeta_1, \bzeta_2-\bk)    
 +
 \hat{M}_{2} (  \bxi_1+\bk, \bxi_2-\bk, \bzeta_1, \bzeta_2)  \\
\nonumber &&\hspace*{0.7in} 
+ 
 \hat{M}_{2} (  \bxi_1+\bk,\bxi_2, \bzeta_1,  \bzeta_2-\bk)     -
2 \hat{M}_{2}(\bxi_1,\bxi_2, \bzeta_1, \bzeta_2)\\
&&\hspace*{0.7in}  -
 \hat{M}_{2} (  \bxi_1,\bxi_2-\bk, \bzeta_1, \bzeta_2-\bk)  
- \hat{M}_{2} (  \bxi_1,\bxi_2+\bk,  \bzeta_1, \bzeta_2-\bk) 
\bigg]d \bk ,
\label{eq:fouriermom0}
\end{eqnarray}
starting from 
$$
\hat{M}_{2}(z=0,\bxi_1,\bxi_2,\bzeta_1,\bzeta_2) = (2\pi r_0 \rho_0)^4 
\exp \Big( - \frac{\rho_0^2}{2} ( |\bxi_1|^2+ |\bxi_2|^2 ) - \frac{r_0^2}{2}( |\bzeta_1|^2+ |\bzeta_2|^2) \Big). 
$$
The solution of this transport equation would give the expression of the variance of the refocused wave.
However, in contrast to the second-order moment of the paraxial Green's function, 
we cannot solve this equation and find a closed-form
expression.
Therefore we address in the next sections a particular regime in which explicit expressions can be obtained.

\section{The Scintillation Regime}
\label{sec:regime}%
In this paper we  address a regime which can be considered
as a particular case of the paraxial white-noise regime:  the scintillation regime. 
The scintillation regime is valid if the transverse correlation length of the Brownian field $B$
(which is the correlation length of the medium)
is smaller than the radius of the TRM and the one of the smoothing kernel.
If the correlation length is our reference length, this means that 
in this regime the covariance function $C^\eps$,
the radius of the TRM $r_0^\eps$, the smoothing kernel radius $\rho_0^\eps$,
and the propagation distance $L^\eps$
are of the form
\begin{equation}
\label{sca:sci}
C^\eps(\bx)= \eps C (\bx) , \quad \quad r_0^\eps = \frac{r_0}{\eps}, \quad \quad \rho_0^\eps = \frac{\rho_0}{\eps},
\quad \quad L^\eps = \frac{L}{\eps} .
\end{equation}
Here $\eps$ is a small dimensionless parameter and we will study the limit $\eps \to 0$.

Let us denote the rescaled function
\begin{equation}
\label{eq:renormhatM2}
\widetilde{M}^\eps (z ,\bxi_1,\bxi_2,\bzeta_1,\bzeta_2) = \hat{M}_2^\eps \Big(\frac{z}{\eps} ,\bxi_1,\bxi_2,\bzeta_1,\bzeta_2 \Big)
 \exp \Big( \frac{i z}{k_0 \eps} (\bxi_2 \cdot \bzeta_2  +   \bxi_1 \cdot \bzeta_1) \Big)
 .
\end{equation}
In the scintillation regime the rescaled function $\widetilde{M}^\eps$  satisfies the equation with fast phases
\begin{eqnarray}
\nonumber
 \frac{\partial \widetilde{M}^\eps}{\partial z} 
&=&
\frac{k_0^2}{4 (2\pi)^2} 
\int \hat{C}(\bk) \bigg[  - 2 \widetilde{M}^\eps(  \bxi_1 ,\bxi_2, \bzeta_1,\bzeta_2)   \\
\nonumber
&&+
\widetilde{M}^\eps (  \bxi_1-\bk,\bxi_2-\bk,  \bzeta_1,\bzeta_2) 
e^{i\frac{z}{\eps k_0} \bk \cdot  (\bzeta_2 + \bzeta_1)} \\
&& \nonumber
+ 
\widetilde{M}^\eps  (  \bxi_1-\bk,\bxi_2,  \bzeta_1,\bzeta_2-\bk) 
e^{i\frac{z}{\eps k_0} \bk \cdot ( \bxi_2 +  \bzeta_1)}\\
\nonumber
&& +\widetilde{M}^\eps (  \bxi_1+\bk,\bxi_2-\bk,  \bzeta_1,\bzeta_2) 
e^{i\frac{z}{\eps k_0} \bk \cdot (\bzeta_2 -   \bzeta_1)}\\
&& \nonumber
+ 
\widetilde{M}^\eps  (  \bxi_1+\bk,\bxi_2,  \bzeta_1,\bzeta_2-\bk) 
e^{i\frac{z}{\eps k_0} \bk \cdot (\bxi_2 -  \bzeta_1)}\\
\nonumber
&&
-
 \widetilde{M}^\eps  (  \bxi_1,\bxi_2-\bk, \bzeta_1, \bzeta_2-\bk) 
 e^{i\frac{z}{\eps k_0} (  \bk \cdot (\bzeta_2+\bxi_2)-|\bk|^2 )} \\
 &&
- \widetilde{M}^\eps  ( \bxi_1,\bxi_2-\bk,  \bzeta_1,\bzeta_2+\bk) 
e^{i\frac{z}{\eps k_0} ( \bk \cdot (\bzeta_2-\bxi_2)+|\bk|^2)}
\bigg] d \bk ,
\label{eq:tildeNeps}
\end{eqnarray}
starting from 
\begin{equation}
\label{eq:initialtildeM2eps}
\widetilde{M}^\eps(z=0,\bxi_1,\bxi_2, \bzeta_1, \bzeta_2 ) = (2\pi)^8 \phi^\eps_{\rho_0} ( \bxi_1 )
\phi^\eps_{\rho_0} ( \bxi_2 )
\phi^\eps_{r_0} ( \bzeta_1 )\phi^\eps_{r_0} ( \bzeta_2 ) ,
\end{equation}
where we have denoted
\begin{equation}
\phi^\eps_{\rho_0}(\bxi) = \frac{\rho_0^2}{2\pi \eps^2} \exp \Big( -\frac{\rho_0^2}{2 \eps^2} |\bxi|^2\Big) ,
\end{equation}
and similarly for $\phi^\eps_{r_0}$.
Note that $\phi^\eps_{\rho_0}$ belongs to $L^1$ and has a $L^1$-norm equal to one,
and that it behaves like a Dirac distribution as $\eps \to 0$.
The following result  shows that $\widetilde{M}^\eps$ exhibits a multi-scale behavior
as $\eps \to 0$, with some components evolving at the scale $\eps$ and 
some components evolving at the scale $1$.
\begin{prop}
\label{prop:sci1}%
The function $\widetilde{M}^\eps(z,\bxi_1,\bxi_2, \bzeta_1,\bzeta_2 ) $ can be expanded as
\begin{eqnarray}
\nonumber
&&\hspace*{-0.2in}
 \widetilde{M}^\eps(z,\bxi_1,\bxi_2,  \bzeta_1,\bzeta_2 )  =
K(z)
\phi^\eps_{\rho_0} ( \bxi_1 )
\phi^\eps_{\rho_0} ( \bxi_2 )
\phi^\eps_{r_0} ( \bzeta_1 )
\phi^\eps_{r_0} ( \bzeta_2 ) \\
\nonumber
&& 
\hspace*{-0.1in}
+
K(z) 
\phi^\eps_{\rho_0} \big( \frac{\bxi_1-\bxi_2}{\sqrt{2}}\big)
\phi^\eps_{r_0} ( \bzeta_1 )
\phi^\eps_{r_0} ( \bzeta_2 )
A\big(z, \frac{\bxi_2+\bxi_1}{2} ,\frac{\bzeta_2 + \bzeta_1}{\eps}\big) \\
\nonumber
&& \hspace*{-0.1in}
+
K(z) 
\phi^\eps_{\rho_0} \big( \frac{\bxi_1+\bxi_2}{\sqrt{2}}\big)
\phi^\eps_{r_0} ( \bzeta_1 )
\phi^\eps_{r_0} ( \bzeta_2 )
A \big(z, \frac{\bxi_2-\bxi_1}{2} ,\frac{\bzeta_2- \bzeta_1}{\eps} \big) \\
\nonumber
&& \hspace*{-0.1in}
+
K(z) 
\phi^\eps_{R_0} \big( \frac{\bxi_1-\bzeta_2}{\sqrt{2}}\big)
\phi^\eps_{r_0} ( \bzeta_1 ) \phi^\eps_{\rho_0} ( \bxi_2 )
A\big(z, \frac{\bzeta_2+\bxi_1}{2} ,\frac{\bxi_2+ \bzeta_1}{\eps}  \big) \\
\nonumber
&& \hspace*{-0.1in}
+
K(z) 
\phi^\eps_{R_0} \big( \frac{\bxi_1+\bzeta_2}{\sqrt{2}}\big)
\phi^\eps_{r_0} ( \bzeta_1 ) \phi^\eps_{\rho_0} ( \bxi_2 )
A\big(z, \frac{\bzeta_2-\bxi_1}{2} ,\frac{\bxi_2- \bzeta_1}{\eps}  \big) \\
\nonumber
&&\hspace*{-0.1in}
+K(z)\phi^\eps_{r_0} ( \bzeta_1 )\phi^\eps_{r_0} (\bzeta_2) 
A \big( z, \frac{\bxi_2+\bxi_1}{2},   \frac{\bzeta_2+ \bzeta_1}{\eps}  \big)
A \big( z, \frac{\bxi_2-\bxi_1}{2},   \frac{\bzeta_2- \bzeta_1}{\eps}  \big) \\
\nonumber
&&\hspace*{-0.1in} +K(z)  \phi^\eps_{r_0} ( \bzeta_1 ) \phi^\eps_{\rho_0} (\bxi_2)
A \big( z, \frac{\bzeta_2+\bxi_1}{2},  \frac{\bxi_2+ \bzeta_1}{\eps}   \big)
A \big( z, \frac{\bzeta_2-\bxi_1}{2},  \frac{\bxi_2- \bzeta_1}{\eps}  \big)
\\
&&\hspace*{-0.1in}
 + R^\eps  (z ,\bxi_1,\bxi_2 ,  \bzeta_1 ,\bzeta_2 )   ,
\label{eq:propsci11}
\end{eqnarray}
where 
\begin{equation}
\frac{1}{R_0^2} = \frac{1}{2} \Big( \frac{1}{r_0^2} + \frac{1}{\rho_0^2} \Big) ,
\end{equation}
the functions $K$ and $A$  are defined by
\begin{eqnarray}
\label{def:K}
K(z) &=& (2\pi)^8 \exp\Big(- \frac{k_0^2}{2} C({\bf 0}) z\Big) , \\
\label{def:A}
A(z,\bxi,\bzeta)  &=&  \frac{1}{2(2\pi)^2}
 \int  \Big[  \exp \Big( \frac{k_0^2}{4} \int_0^z C\big( \bx + \frac{\bzeta}{k_0} z' \big) dz' \Big) -1\Big]
   \exp \big( -i \bxi\cdot \bx  \big)
 d\bx  ,   \hspace*{0.1in}
\end{eqnarray}
and the function $R^\eps $ satisfies
\begin{equation}
\sup_{z \in [0,Z]} \| R^\eps (z,\cdot,\cdot,\cdot,\cdot) \|_{L^1(\RR^2\times \RR^2\times \RR^2\times \RR^2)} 
\stackrel{\eps \to 0}{\longrightarrow}  0  ,
\end{equation}
for any $Z>0$.
\end{prop}

This result is an extension of Proposition 6.1 in \cite{garniers4} in which the case $r_0=\rho_0$ is addressed.
It shows that, if we deal with an integral of $\widetilde{M}^\eps$ against a bounded function,
then we can replace $\widetilde{M}^\eps$ by the right-hand side of (\ref{eq:propsci11}) without the $R^\eps$ term
up to a negligible error when $\eps$ is small. This substitution will allow us to get explicit 
and quantitative results for the refocusing properties of the SLM scheme.

\section{Refocusing at the Original Source Location}
\label{sec:refocus}%
In the scintillation regime the mean refocused wave is given by (see (\ref{eq:mean1}) 
and take $C \to \eps C$, $r_0 \to r_0/\eps$, $\rho_0 \to \rho_0/\eps$, $\by \to \by/\eps$, $L\to L/\eps$):
 \begin{eqnarray}
 \nonumber
 \EE \Big[ \hat{u}_{\rm tr}^\eps( \frac{\by}{\eps}+\bx ; \frac{\by}{\eps} )\Big] 
 &=&
 \frac{r_0^2}{4 \pi \eps^2} 
 \int \exp \Big( i \bxi \cdot \big(\frac{\by}{\eps}+\frac{\bx}{2}\big) 
  - \frac{r_0^2 |\bxi|^2}{4 \eps^2} -\frac{ |\eps \bx - \bxi \frac{L}{k_0 } |^2}{4 \rho_0^2} \\
 \nonumber
 &&\hspace*{0.7in}
 +\frac{k_0^2 \eps}{4} \int_0^{L/\eps} C ( \bx -\bxi \frac{z}{k_0} ) -C({\bf 0}) dz \Big) d\bxi \\
 \nonumber
&=&
 \frac{r_0^2}{4 \pi} 
 \int \exp \Big( i \bxi \cdot \big( \by +\frac{\eps \bx}{2}\big) 
  - \frac{r_0^2 |\bxi|^2}{4} -\eps^2 \frac{ |\bx - \bxi \frac{L}{k_0 } |^2}{4 \rho_0^2} \\
 &&\hspace*{0.7in}
 +\frac{k_0^2}{4} \int_0^{L} C ( \bx -\bxi \frac{z}{k_0} ) -C({\bf 0}) dz \Big) d\bxi .
  \end{eqnarray}
  
 The following result is then straightforward.
 
 \begin{proposition}
 \label{prop:momtr1}%
The mean refocused wave converges as $\eps \to 0$:
 \begin{equation}
 \EE \Big[ \hat{u}_{\rm tr}^\eps\big( \frac{\by}{\eps}+\bx ; \frac{\by}{\eps} \big)\Big] 
 \stackrel{\eps \to 0}{\longrightarrow} 
 \frac{r_0^2}{4 \pi} 
 \int \exp \Big( i \bxi \cdot  \by 
  - \frac{r_0^2 |\bxi|^2}{4}
 +\frac{k_0^2}{4} \int_0^{L} C ( \bx -\bxi \frac{z}{k_0} ) -  C({\bf 0}) dz \Big) d\bxi .
 \label{eq:mainpeak}
  \end{equation}
  \end{proposition}

We can observe  that, if the number of elements of the TRM is
large, i.e. $\rho_0 \ll R_{\rm m}$, then $r_0 \simeq R_{\rm m}$ and the profile of the mean refocused wave
does not depend anymore on the number and size of the elements, but only on the radius of the TRM.
The expression (\ref{eq:mainpeak}) shows that, as a function of the offset $\bx$, the mean refocused wave 
has the form of a peak centered at $\bx={\bf 0}$. This peak has maximal amplitude
 \begin{equation}
{\cal U}_{{\rm b},\by} =  \lim_{\eps \to 0}
  \EE \Big[ \hat{u}_{\rm tr}^\eps\big( \frac{\by}{\eps}  ; \frac{\by}{\eps} \big)\Big] 
=
 \frac{r_0^2}{4 \pi} 
 \int \exp \Big( i \bxi \cdot  \by 
  - \frac{r_0^2 |\bxi|^2}{4}
 +\frac{k_0^2}{4} \int_0^{L} C (\bxi \frac{z}{k_0} ) -  C({\bf 0}) dz \Big) d\bxi  ,
 \label{eq:mainpeakp}
  \end{equation}
and it raises above the  constant background
 \begin{equation}
{\cal U}_{{\rm p},\by} =  \lim_{|\bx|\to \infty}
\lim_{\eps \to 0} \EE \Big[ \hat{u}_{\rm tr}^\eps\big( \frac{\by}{\eps}+\bx ; \frac{\by}{\eps} \big)\Big] 
=\exp \Big( - \frac{|\by|^2}{r_0^2}
  - \frac{k_0^2   C({\bf 0})  L}{4}   \Big)  .
 \label{eq:mainpeakb}
  \end{equation}


By (\ref{eq:secondmom1b}), in the scintillation regime, the second moment 
of the refocused wave is
\begin{eqnarray}
\nonumber
 &&\EE \Big[  \hat{u}_{\rm tr}^\eps \big( \frac{\by}{\eps}+\bx +\frac{\bh}{2}; \frac{\by}{\eps} \big) 
 \overline{\hat{u}_{\rm tr}^\eps \big( \frac{\by}{\eps}+\bx-\frac{\bh}{2} ; \frac{\by}{\eps} \big)}  \Big] \\
\nonumber && =
 M_2^\eps \big( \frac{L}{\eps}, \br_1= \bx+ \frac{2 \by}{\eps}, \br_2= \frac{\bh}{2} , \bq_1=\frac{\bh}{2}, \bq_2 = \bx \big) \\
\nonumber&& =
 \frac{1}{(2\pi)^8} \iint 
 \widetilde{M}^\eps (L, \bxi_1,\bxi_2,\bzeta_1,\bzeta_2)
 \exp  \Big( i (\bx+ \frac{2 \by}{\eps})\cdot \bzeta_1+i  \frac{\bh}{2} \cdot ( \bzeta_2+\bxi_1) + i\bx \cdot \bxi_2 \Big)  \\
&& \hspace*{0.8in} \times  \exp \Big( - \frac{i L}{k_0 \eps} (\bxi_2 \cdot \bzeta_2  +   \bxi_1 \cdot \bzeta_1) \Big)
 d\bxi_1d\bxi_2 d\bzeta_1d\bzeta_2  .
\end{eqnarray}
By Proposition \ref{prop:sci1}, we find that it converges as $\eps \to 0$ to
\begin{eqnarray}
\nonumber
 \EE \Big[  \hat{u}_{\rm tr}^\eps \big( \frac{\by}{\eps}+\bx +\frac{\bh}{2}; \frac{\by}{\eps} \big) 
 \overline{\hat{u}_{\rm tr}^\eps \big( \frac{\by}{\eps}+\bx-\frac{\bh}{2} ; \frac{\by}{\eps} \big)}  \Big] 
 & \stackrel{\eps \to 0}{\longrightarrow} &
\frac{K}{(2\pi)^8} \int 
\phi^1_{r_0} ( \bzeta_1 )
e^{ 2 i \bzeta_1 \cdot \by } d\bzeta_1 \\
\nonumber
&& 
\hspace*{-2.7in}
+\frac{2 K}{(2\pi)^8} \iint 
\phi^1_{r_0} ( \bzeta_1 )
\phi^1_{r_0} ( \bzeta_2 )
\big[ 
A\big( \bxi_2 , \bzeta_2 + \bzeta_1 \big) 
e^{ 2 i \bzeta_1 \cdot \by -i \frac{L}{k_0} (\bzeta_2 +\bzeta_1) \cdot \bxi_2 +i \bxi_2 \cdot (\bx+\frac{\bh}{2})} 
\\
\nonumber
&& \hspace*{-1.5in}
+
A\big(\bxi_2 , \bzeta_2 - \bzeta_1 \big) 
e^{ 2 i \bzeta_1 \cdot \by -i \frac{L}{k_0} (\bzeta_2 -\bzeta_1) \cdot \bxi_2 +i \bxi_2 \cdot (\bx-\frac{\bh}{2})} 
\big] d\bzeta_1 d\bzeta_2 d\bxi_2 \\
\nonumber
&& \hspace*{-2.7in}
+\frac{2 K}{(2\pi)^8} \iint 
\phi^1_{r_0} ( \bzeta_1 )
\phi^1_{\rho_0} ( \bxi_2 )
\big[ 
A\big(\bzeta_2 , \bxi_2 + \bzeta_1 \big) 
e^{ 2 i \bzeta_1 \cdot \by -i \frac{L}{k_0} (\bxi_2 +\bzeta_1) \cdot \bzeta_2 + i \bzeta_2\cdot \bh } 
\\
\nonumber
&& \hspace*{-1.5in}
+
A\big( \bzeta_2 , \bxi_2 - \bzeta_1 \big) 
e^{ 2 i \bzeta_1 \cdot \by -i \frac{L}{k_0} (\bxi_2 -\bzeta_1) \cdot \bzeta_2 } 
\big] d\bzeta_1 d\bzeta_2 d\bxi_2
 \\
\nonumber
&&\hspace*{-2.7in}
+\frac{K}{(2\pi)^8} \iint 
\phi^1_{r_0} ( \bzeta_1 )
\phi^1_{r_0} ( \bzeta_2 )
A\big( \frac{\bxi_2 +\bxi_1}{2}, \bzeta_2 + \bzeta_1 \big) 
A\big( \frac{\bxi_2 -\bxi_1}{2} , \bzeta_2 - \bzeta_1 \big) \\
\nonumber
&& \hspace*{-1.5in}\times
e^{ 2 i \bzeta_1 \cdot \by +i \bxi_2\cdot \bx -i \frac{L}{k_0} (\bxi_1 \cdot \bzeta_1+\bxi_2 \cdot \bzeta_2 )
+i \bxi_1 \cdot \frac{\bh}{2}  } 
d\bzeta_1 d\bzeta_2d \bxi_1 d\bxi_2
\\
\nonumber
&&\hspace*{-2.7in}
+\frac{K}{(2\pi)^8} \iint 
\phi^1_{r_0} ( \bzeta_1 )
\phi^1_{\rho_0} ( \bxi_2 )
A\big( \frac{\bzeta_2 +\bxi_1}{2}, \bxi_2 + \bzeta_1 \big) 
A\big( \frac{\bzeta_2 -\bxi_1}{2} , \bxi_2 - \bzeta_1 \big) \\
&& \hspace*{-1.5in}\times
e^{ 2 i \bzeta_1 \cdot \by  -i \frac{L}{k_0} (\bxi_1 \cdot \bzeta_1+\bxi_2 \cdot \bzeta_2 )
+ i (\bzeta_2+\bxi_1) \cdot \frac{\bh}{2}} 
d\bzeta_1 d\bzeta_2d \bxi_1 d\bxi_2.
\end{eqnarray}
Using the explicit expressions (\ref{def:K}) and (\ref{def:A}) for  $K$ and $A$, we get
\begin{eqnarray}
\nonumber
\EE \Big[  \hat{u}_{\rm tr}^\eps \big( \frac{\by}{\eps}+\bx +\frac{\bh}{2}; \frac{\by}{\eps} \big) 
 \overline{\hat{u}_{\rm tr}^\eps \big( \frac{\by}{\eps}+\bx-\frac{\bh}{2} ; \frac{\by}{\eps} \big)}  \Big] 
 & \stackrel{\eps \to 0}{\longrightarrow} &
- \exp \Big( -\frac{k_0^2 C({\bf 0}) L}{2} - \frac{2 |\by|^2}{r_0^2} \Big)\\
\nonumber
&& \hspace*{-2.7in}
+
\Big( \frac{r_0^2}{4 \pi} \Big)^2 
\Big[ \int \exp \Big( - \frac{r_0^2}{4} |\balpha|^2 + i \balpha \cdot \by + \frac{k_0^2}{4} \int_0^L C( \bx +\frac{\bh}{2} - \balpha \frac{z}{k_0} ) -C({\bf 0}) dz \Big) d\balpha \Big]\\
\nonumber
&& \hspace*{-2.3in} \times
\Big[ \int \exp \Big( - \frac{r_0^2}{4} |\balpha|^2 - i \balpha \cdot \by + \frac{k_0^2}{4} \int_0^L C( \bx -\frac{\bh}{2} - \balpha \frac{z}{k_0} ) -C({\bf 0}) dz \Big) d\balpha \Big]\\
\nonumber
&& \hspace*{-2.7in}
+
\Big( \frac{r_0\rho_0}{4 \pi} \Big)^2 
  \int \exp \Big( - \frac{r_0^2+\rho_0^2}{8} (|\balpha|^2+|\bbeta|^2) + \frac{r_0^2-\rho_0^2}{4} \balpha \cdot \bbeta 
+  i (\balpha-\bbeta) \cdot \by \Big) \\
&&  \hspace*{-1.9in}\times \exp \Big( - \frac{k_0^2 C({\bf 0}) L}{2} + \frac{k_0^2}{4} \int_0^L C(\bh- \balpha \frac{z}{k_0} ) + C(\bbeta \frac{z}{k_0} )  dz \Big) d\balpha d\bbeta 
   .
\end{eqnarray}
As a consequence we can now describe the covariance function of the refocused field at offsets $\bx +{\bh}/{2}$
and $\bx -{\bh}/{2}$ (relatively to the original source location $\by/\eps$) defined by:
\begin{eqnarray*}
&&{\rm Cov} \Big(  \hat{u}_{\rm tr}^\eps \big( \frac{\by}{\eps}+\bx +\frac{\bh}{2}; \frac{\by}{\eps} \big) ,
 \hat{u}_{\rm tr}^\eps \big( \frac{\by}{\eps}+\bx-\frac{\bh}{2} ; \frac{\by}{\eps} \big)  \Big) \\
&& =
 \EE \Big[  \hat{u}_{\rm tr}^\eps \big( \frac{\by}{\eps}+\bx +\frac{\bh}{2}; \frac{\by}{\eps} \big) 
 \overline{\hat{u}_{\rm tr}^\eps \big( \frac{\by}{\eps}+\bx-\frac{\bh}{2} ; \frac{\by}{\eps} \big)}  \Big] \\
&&\quad -
\EE \Big[  \hat{u}_{\rm tr}^\eps \big( \frac{\by}{\eps}+\bx +\frac{\bh}{2}; \frac{\by}{\eps} \big) \Big]
\EE \Big[ \overline{\hat{u}_{\rm tr}^\eps \big( \frac{\by}{\eps}+\bx-\frac{\bh}{2} ; \frac{\by}{\eps} \big)}  \Big]   .
\end{eqnarray*}

\begin{proposition}
 \label{prop:momtr2}%
The covariance function of the refocused field satisfies
\begin{eqnarray}
\nonumber
&& {\rm Cov} \Big(  \hat{u}_{\rm tr}^\eps \big( \frac{\by}{\eps}+\bx +\frac{\bh}{2}; \frac{\by}{\eps} \big) ,
 \hat{u}_{\rm tr}^\eps \big( \frac{\by}{\eps}+\bx-\frac{\bh}{2} ; \frac{\by}{\eps} \big)  \Big) \\
\nonumber
  && \stackrel{\eps \to 0}{\longrightarrow} 
 \Big( \frac{r_0\rho_0}{4 \pi} \Big)^2 
  \iint \exp \Big( - \frac{r_0^2+\rho_0^2}{8} (|\balpha|^2+|\bbeta|^2) + \frac{r_0^2-\rho_0^2}{4} \balpha \cdot \bbeta 
+  i (\balpha-\bbeta) \cdot \by  \Big) \\
&& \quad\quad \times \exp \Big( - \frac{k_0^2 C({\bf 0}) L}{2}\Big)
\Big[ \exp \Big(  \frac{k_0^2}{4} \int_0^L C(\bh-\balpha \frac{z}{k_0} ) + C(\bbeta \frac{z}{k_0} )   dz \Big)-1 \Big] d\balpha d\bbeta . \hspace*{0.2in}
\end{eqnarray}
\end{proposition}
Note that the covariance does not depend on the central off-set $\bx$. 
This means in particular that the variance of  the refocused field is constant 
in the neighborhood of the original source location.

The refocused wave $\bx\to  \hat{u}_{\rm tr}^\eps \big( {\by}/{\eps} +\bx ; {\by}/{\eps} \big)$
therefore consists of a main peak centered at $\bx={\bf 0}$, of the form
(\ref{eq:mainpeak}), with peak intensity (square difference between the maximal amplitude and the background amplitude)
\begin{eqnarray}
\nonumber
{\cal I}_{{\rm p},\by} &=&
\big| {\cal U}_{{\rm b},\by}-{\cal U}_{{\rm p},\by}\big|^2  \\
\nonumber
&=&   
\Big( \frac{r_0^2}{4 \pi} \Big)^2 \Big| \int \exp \Big( i \bxi \cdot  \by 
  - \frac{r_0^2 |\bxi|^2}{4} \Big) \\
  &&   \times \exp \Big( - \frac{k_0^2 C({\bf 0}) L}{4} \Big) 
 \Big[ \exp 
 \Big( \frac{k_0^2}{4} \int_0^{L} C ( \bxi \frac{z}{k_0} )   dz   \Big)  -1\Big]
 d\bxi \Big|^2    ,
\label{eq:Ip1}
\end{eqnarray}
over a zero-mean fluctuating background with intensity
\begin{eqnarray}
\nonumber
{\cal I}_{{\rm b},\by} &=& \lim_{\eps \to 0} {\rm Var} \Big(   \hat{u}_{\rm tr}^\eps \big( \frac{\by}{\eps}+\bx ; \frac{\by}{\eps} \big)   \Big) =
 \lim_{\eps \to 0}  \EE \Big[   \big| \hat{u}_{\rm tr}^\eps \big( \frac{\by}{\eps} +\bx ; \frac{\by}{\eps} \big) \big|^2 \Big]  
 -  \Big| \EE \Big[   \hat{u}_{\rm tr}^\eps \big( \frac{\by}{\eps} +\bx; \frac{\by}{\eps} \big) \Big]  \Big|^2 \\
\nonumber
&=& \Big( \frac{r_0\rho_0}{4 \pi} \Big)^2 
  \iint \exp \Big( - \frac{r_0^2+\rho_0^2}{8} (|\balpha|^2+|\bbeta|^2) + \frac{r_0^2-\rho_0^2}{4} \balpha \cdot \bbeta 
+  i (\balpha-\bbeta) \cdot \by  \Big) \\
&& \quad \times \exp \Big( - \frac{k_0^2 C({\bf 0}) L}{2}\Big)
\Big[ \exp \Big(  \frac{k_0^2}{4} \int_0^L C(\balpha \frac{z}{k_0} ) + C(\bbeta \frac{z}{k_0} )   dz \Big)-1 \Big] d\balpha d\bbeta .
\label{eq:Ib1}
\hspace*{0.2in}
\end{eqnarray} 
 Let us   assume that the covariance function $C$ is isotropic and at least twice differentiable at zero and
write it in the form
\begin{equation}
\label{eq:formcov}
C(\bx) = \sigma^2 l_c \widetilde{C} \Big( \frac{|\bx|}{l_c} \Big) ,
\end{equation}
with $ \widetilde{C}(0)=1$, $\widetilde{C}'(0)=0$, and $\widetilde{C}''(0)=-1$. In this framework $C({\bf 0})=\sigma^2 l_c$ and 
the correlation radius of the medium is $l_c$.
When the original source is at $\by={\bf 0}$, then the peak intensity is
\begin{eqnarray}
\nonumber
{\cal I}_{{\rm p},{\bf 0}} &=& \exp \Big( -\frac{\sigma^2 k_0^2 l_c L}{2} \Big) 
 \Big| \int_0^\infty a \exp \Big(  
  - \frac{a^2}{2} \Big) \\
  && \times \Big[ \exp \Big( \frac{\sigma^2 k_0^2 l_c L}{4} \int_0^{1} \widetilde{C} \big(a \frac{\sqrt{2} L}{k_0l_c r_0} s \big) ds \Big) -1 \Big]  da \Big|^2 ,
\label{eq:Ip2}
\hspace*{0.25in}
\end{eqnarray}
and the background intensity is
\begin{eqnarray}
\nonumber
&& {\cal I}_{{\rm b},{\bf 0}} = \Big( \frac{2 r_0 \rho_0}{r_0^2 +\rho_0^2} \Big)^2 
 \exp \Big( -\frac{\sigma^2 k_0^2 l_c L}{2} \Big) 
\int_0^\infty \int_0^\infty ab \exp \Big(  
  - \frac{a^2+b^2}{2} \Big) I_0 \Big( \frac{ r_0^2 - \rho_0^2}{r_0^2 +\rho_0^2} ab \Big) \\
&& 
\times \Big[ 
\exp \Big( \frac{\sigma^2 k_0^2 l_c L}{4} \int_0^{1} \widetilde{C} \big(a \frac{2L}{k_0l_c \sqrt{r_0^2+\rho_0^2}}s \big) 
+ \widetilde{C} \big(b \frac{2L}{k_0l_c \sqrt{r_0^2+\rho_0^2}}s \big) ds \Big)- 1 \Big]  da db , \hspace*{0.3in}
\label{eq:Ib2}
\end{eqnarray}
where $I_0$ is the modified  Bessel function of the first kind and order zero.
As a consequence we can now describe the SNR of the refocused field.

\begin{proposition}
If the covariance function $C$ is of the form (\ref{eq:formcov}),
then the contrast or signal-to-noise ratio of the refocused field defined by
\begin{equation}
\label{def:SNRTR}
{\rm SNR}_{\rm tr} =
 \lim_{\eps \to 0} 
\frac{\big| \EE \big[   \hat{u}_{\rm tr}^\eps({\bf 0}  ; {\bf 0}  ) \big]    - \lim_{|\bx'|\to \infty}  \EE \big[   \hat{u}_{\rm tr}^\eps(\bx' ; {\bf 0}  ) \big]   \big|^2 }
{ {\rm Var} \big( \hat{u}_{\rm tr}^\eps( \bx  ; {\bf 0}  )  \big)  }
,
\end{equation}
does not depend on $\bx$ and it is equal to 
\begin{equation}
{\rm SNR}_{\rm tr} = \frac{{\cal I}_{{\rm p},{\bf 0}}  }{{\cal I}_{{\rm b},{\bf 0}} }  ,
\end{equation}
where ${\cal I}_{{\rm p},{\bf 0}} $
 is given by (\ref{eq:Ip2})
and ${\cal I}_{{\rm b},{\bf 0}} $
 is given by (\ref{eq:Ib2}).
 \end{proposition}

In order to make the discussion more explicit, 
we assume in the following that scattering is strong in the sense that the propagation distance
 is larger than the scattering mean free path $\sigma^2 k_0^2 l_c L \gg 1$.
The scattering mean free path $L_{\rm sca} = 4/(\sigma^2 k_0^2 l_c )$ determines
the exponential decay rate of the mean amplitude of a wave propagating through the random medium,
which decays as $\exp (- k_0^2 C({\bf 0}) L/8) = \exp (- \sigma^2 k_0^2 l_c L/8) = \exp (-L/(2L_{\rm sca})) $ 
with $L_{\rm sca} = 4/(\sigma^2 k_0^2 l_c )$,  as can be seen from the It\^o's form  (\ref{eq:modelito})).
Using  (\ref{eq:mainpeak})  with  $\by={\bf 0}$, we then find that
the main peak is  a Gaussian peak centered at ${\bf 0}$,
\begin{equation}
 \EE \big[ \hat{u}_{\rm tr}^\eps \big( \bx ; {\bf 0}\big)\big] 
 \stackrel{\eps \to 0}{\longrightarrow} 
 \frac{1}{  1 +  \frac{\sigma^2  L^3}{6 r_0^2 l_c} } 
\exp \Big( - \frac{1 +  \frac{\sigma^2 L^3}{24 r_0^2 l_c} }{1 +  \frac{\sigma^2  L^3}{6 r_0^2 l_c} }
 \frac{\sigma^2  k_0^2 L}{8 l_c} |\bx|^2 \Big)  ,
\end{equation}
which is independent of $\rho_0$.
 The width of the peak is $R_{\rm tr} $ given by
\begin{equation}
R_{\rm tr}^2 = \frac{4 l_c}{\sigma^2  k_0^2  L}  \frac{1 + \frac{\sigma^2  L^3}{6 r_0^2 l_c}}{1 + \frac{\sigma^2  L^3}{24 r_0^2 l_c}} ,
\label{def:widthtr}
\end{equation}
and the peak intensity is
\begin{equation}
{\cal I}_{{\rm p},{\bf 0}}  =  \frac{1}{  \big( 1 +  \frac{\sigma^2  L^3}{6 r_0^2l_c}\big)^2  } .
\end{equation}
From (\ref{eq:Ib1}) the background intensity is
\begin{equation}
{\cal I}_{{\rm b},{\bf 0}} =  \frac{1}{  \big( 1 +  \frac{\sigma^2  L^3}{6 r_0^2l_c}\big) \big( 1 +  \frac{\sigma^2  L^3}{6 \rho_0^2 l_c}\big)  }  .
\end{equation}
As a result the signal-to-noise ratio is given by
\begin{equation}
\label{eq:snr1}
{\rm SNR}_{\rm tr} = \frac{1 + \frac{\sigma^2  L^3}{6 \rho_0^2 l_c}}{1 + \frac{\sigma^2  L^3}{6 r_0^2 l_c}}  ,
\end{equation}
and it depends on $\rho_0$. 

We may identify three situations (remember we always have $r_0 \geq \rho_0$, and we are mainly interested
in the case $r_0 \gg \rho_0 $):
\begin{equation}
{\rm SNR}_{\rm tr} \approx
\left\{
\begin{array}{ll}
\displaystyle
1 & 
\displaystyle
\mbox{ if } \rho_0^2  > \frac{\sigma^2  L^3}{6 l_c} ,\\
\displaystyle
\frac{\sigma^2  L^3}{6 \rho_0^2 l_c} & 
\displaystyle
\mbox{ if } \rho_0^2 < \frac{\sigma^2  L^3}{6l_c}< r_0^2 ,\\
\displaystyle
\frac{r_0^2}{\rho_0^2} &
\displaystyle
 \mbox{ if } r_0^2 < \frac{\sigma^2  L^3}{6 l_c} .
\end{array} \right.
\end{equation}
This shows that stability increases when the number of elements of the TRM increases,
or equivalently when the radius $\rho_0$ decreases.
As a function of $r_0$ the signal-to-noise ratio is maximal 
when $ r_0^2 < \sigma^2  L^3/(6l_c)$, and then it
 is given by the number of elements of the TRM $r_0^2/\rho_0^2$.
 We can explain the physical origin of the condition $ r_0^2 < \sigma^2  L^3/(6l_c)$ as follows.
 The field generated by the original source at the target point $({\bf 0},L)$, transmitted through the medium,
 and recorded by the TRM, has the form of a diffuse beam with radius of the order of $\sigma^2 L^3/l_c$ \cite{garniers1}.
 If the radius of the TRM is larger, then the elements of the TRM outside the support of the diffuse beam do
 not record anything and cannot participate in the time-reversal process. They are not used at all,
 that is why the SNR cannot reach its maximal value given by the number of elements of the TRM $r_0^2/\rho_0^2$.

\section{Focusing at a Prescribed Target Point}
\label{sec:prescribed}%
\begin{figure}
\begin{center}
\begin{picture}(300,80)
\put(0,0){\includegraphics[width=6.1cm]{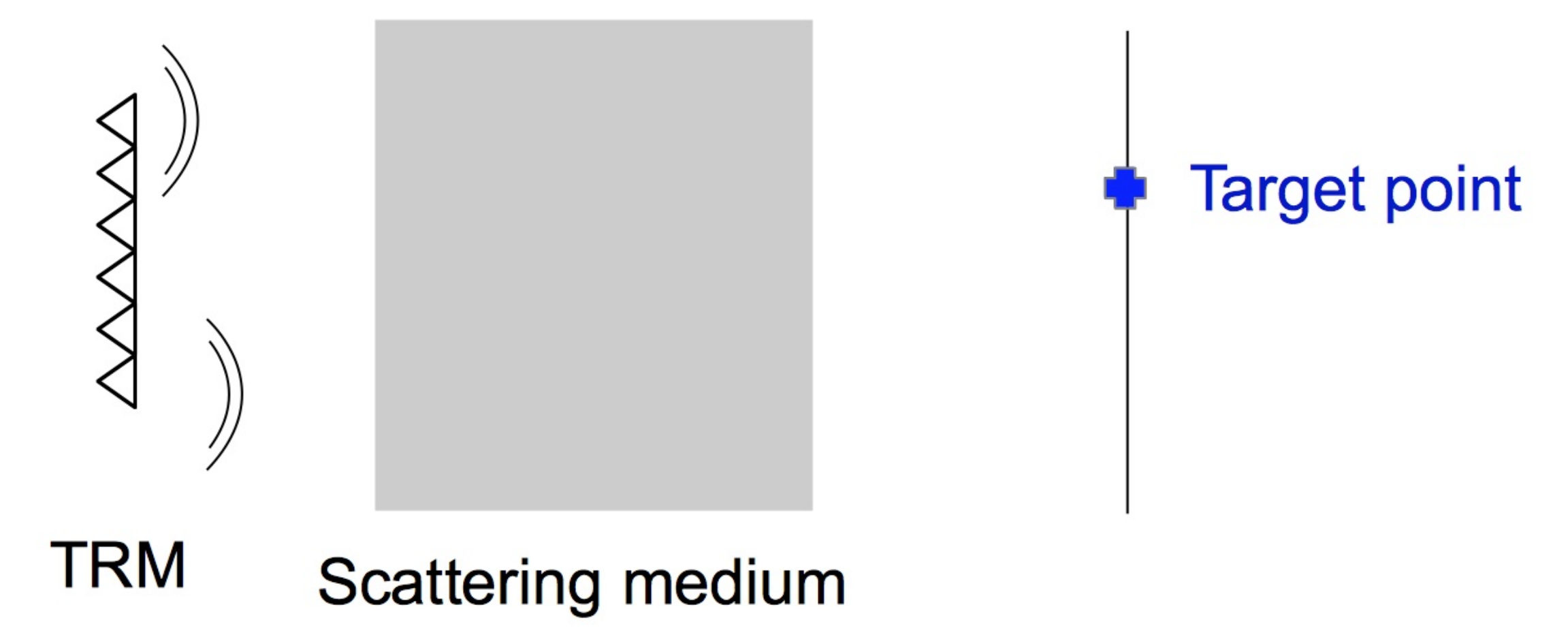}}
\put(170,-4){ \includegraphics[width=4.0cm]{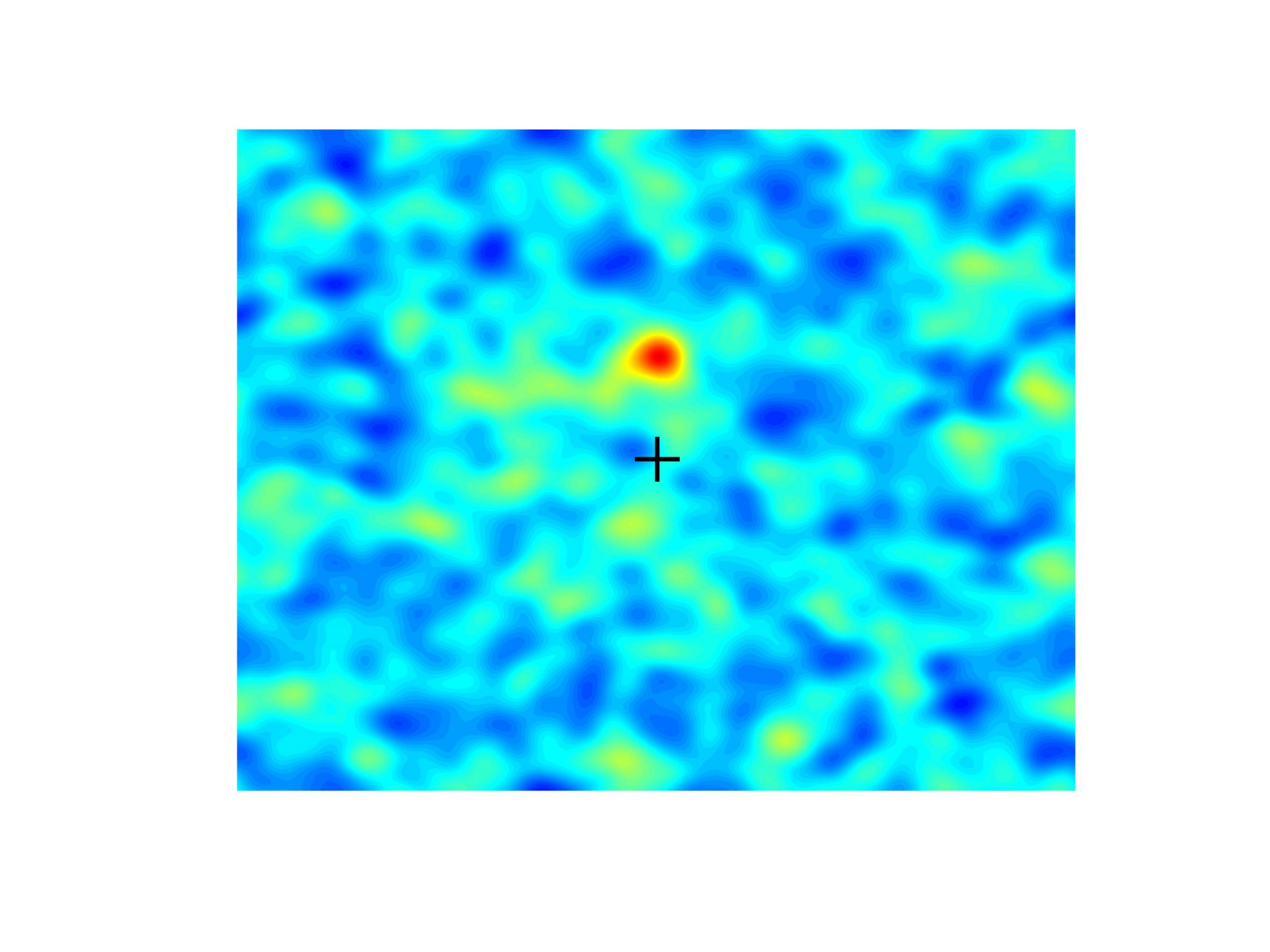}}
\end{picture}
\end{center}
\caption{Focusing on a prescribed point in the neighborhood of the original source location.
The TRM emits the complex-conjugated recorded field with an additional linear phase (the cross in the right image
stands for the original source location).
\label{fig:0c} 
}
\end{figure}

It was noticed in the experiments that wave focusing could be achieved with the optimized 
phases of the SLM in the neighborhood of the target point and not only on the target point \cite{katz12,popoff10}.
In the time-reversal context, this amounts to say that, for a given original point source at $(\by ,L)$,
it is possible to manipulate the field emitted by the TRM to focus on a point 
prescribed by the user near $(\by , L)$.
In this section we introduce and discuss this manipulation 
and we explain and quantify the refocusing property on a prescribed point (see Figure \ref{fig:0c}).

In order to shift the focal spot of the refocused wave the idea is to impose an additional linear phase 
at the TRM, so that the field in the plane $z=L$ is instead of (\ref{eq:trfield1}):
\begin{equation}
\label{eq:trfield2}
\hat{u}_{\rm tr}^{\bb}( \bx ;\by ) =
 \int \exp \Big( - \frac{|\bx_m|^2}{ R_{\rm m}^2} + \frac{i \bb \cdot \bx_m}{  R_{\rm m}^2} \Big) \hat{u}_{\rm em} (\bx,\bx_m)
\overline{ \hat{u}_{\rm rec}(\bx_m;\by) }  
  d\bx_m ,
\end{equation}
where $\hat{u}_{\rm rec}$ has been obtained with a point source at $(\by,L)$ and is given by (\ref{def:urec}).
The vector $\bb$ determines the linear phase and we will see below how to choose it to focus on a prescribed point.
The  time-reversed field can therefore be expressed as
\begin{eqnarray}
\nonumber
\hat{u}_{\rm tr}^{\bb}( \bx ;\by ) = 4 k_0^2 C_0 \exp \Big( - \frac{\rho_0^2}{4 r_0^2(r_0^2-\rho_0^2)} |\bb|^2 \Big)
\iint 
\exp \Big( - \frac{|\bx'|^2}{r_0^2} - \frac{|\by'|^2}{4 \rho_0^2} + i \frac{\bb \cdot \bx'}{r_0^2} \Big) \\
\times \hat{G}\big(L,\bx'+\frac{\by'}{2},\bx\big) 
\overline{\hat{G}\big(L,\bx'-\frac{\by'}{2},\by\big) }
  d\bx' d\by'  , \hspace*{0.25in}
\end{eqnarray}
with $C_0$ and $r_0$ defined by (\ref{def:r0}). We will take $C_0=1$ as in the previous sections.

We proceed as in the previous section.
We consider the scintillation regime $C \to \eps C$, $r_0 \to r_0/\eps$, $\rho_0 \to \rho_0/\eps$, 
$\by \to \by/\eps$, $L \to L/\eps$, $\bb \to \bb /\eps$.
The mean refocused wave at offset $\bx$ (relatively to the original source location $\by/\eps$)  is then of the form
 \begin{eqnarray}
 \nonumber
 \EE \Big[ \hat{u}_{\rm tr}^{\bb,\eps}\big( \frac{\by}{\eps}+\bx ; \frac{\by}{\eps} \big)\Big] 
 \stackrel{\eps \to 0}{\longrightarrow} 
 \frac{r_0^2}{4 \pi}  \exp \Big( - \frac{\rho_0^2}{4 r_0^2(r_0^2-\rho_0^2)} |\bb|^2 \Big) 
  \int \exp \Big( i \bxi \cdot  \by 
  - \frac{r_0^2 }{4}\big|\bxi - \frac{\bb}{r_0^2} \big|^2 \Big)\\
   \times 
\exp \Big( \frac{k_0^2}{4} \int_0^{L} C \big( \bx -\bxi \frac{z}{k_0} \big) -  C({\bf 0}) dz \Big) d\bxi .\hspace*{0.25in}
 \label{eq:mainpeak2}
  \end{eqnarray}
The variance function of the refocused field at offset $\bx$ satisfies
\begin{eqnarray}
\nonumber
&& {\rm Var} \Big(  \hat{u}_{\rm tr}^{\bb,\eps} \big( \frac{\by}{\eps}+\bx  ; \frac{\by}{\eps} \big)    \Big)
=   \EE \Big[   \big| \hat{u}_{\rm tr}^{\bb,\eps} \big( \frac{\by}{\eps} +\bx ; \frac{\by}{\eps} \big) \big|^2 \Big]  
 -  \Big| \EE \Big[   \hat{u}_{\rm tr}^{\bb,\eps} \big( \frac{\by}{\eps} +\bx; \frac{\by}{\eps} \big) \Big]  \Big|^2 
 \\
\nonumber
  && \stackrel{\eps \to 0}{\longrightarrow} 
 \Big( \frac{r_0\rho_0}{4 \pi} \Big)^2  \exp \Big( - \frac{\rho_0^2}{2 r_0^2(r_0^2-\rho_0^2)} |\bb|^2- \frac{k_0^2 C({\bf 0}) L}{2} \Big)\\
 \nonumber
&& \quad\quad\times  \iint \exp \Big( - \frac{r_0^2+\rho_0^2}{8} (|\balpha|^2+|\bbeta|^2) + \frac{r_0^2-\rho_0^2}{4} \balpha \cdot \bbeta 
+  i (\balpha-\bbeta) \cdot \by  \Big) \\
&& \quad\quad\quad\quad \times 
\Big[ \exp \Big(  \frac{k_0^2}{4} \int_0^L C(\balpha \frac{z}{k_0} ) + C(\bbeta \frac{z}{k_0} )   dz \Big)-1 \Big] d\balpha d\bbeta  .
\end{eqnarray}

Assume that the covariance function $C(\bx)$ has the form (\ref{eq:formcov}).
When scattering is strong in the sense that the propagation distance
 is larger than the scattering mean free path $\sigma^2  k_0^2 l_c L \gg 1$, 
the main peak is of the form of a damped and shifted Gaussian peak
\begin{equation}
 \EE \big[ \hat{u}_{\rm tr}^{\bb,\eps} ( \bx ; {\bf 0}) \big] 
\stackrel{\eps \to 0}{\longrightarrow} 
 \frac{1}{  1 +  \frac{\sigma^2  L^3}{6 r_0^2 l_c} } \exp \Big( - \frac{|\bb|^2}{4 r_0^2} \big( 
 \frac{\frac{\sigma^2  L^3}{24l_c}}{r_0^2 + \frac{\sigma^2  L^3}{24 l_c}} + \frac{\rho_0^2}{r_0^2-\rho_0^2} \big)\Big) 
 \exp \Big( - \frac{|\bx-\bx^\bb|^2}{2 R_{\rm tr}^2} \Big) ,
\end{equation}
whose center  is at 
\begin{equation}
\label{eq:phase2a}
\bx^\bb = \alpha_L \bb, \quad \quad \alpha_L = \frac{L}{2 k_0 r_0^2 (1+\frac{\sigma^2  L^3}{24 r_0^2 l_c})} ,
\end{equation}
and the width $R_{\rm tr}$ is given by (\ref{def:widthtr}).
This shows that the linear phase $\bb$ in (\ref{eq:trfield2}) generates a shift in the focal spot,
that is deterministic, proportional to the linear phase, and fully predictible.
If $\frac{\sigma^2  L^3}{24 l_c} \ll r_0^2$, then $\alpha_L=L/(2k_0 r_0^2)$ and therefore
the linear phase that one needs to impose to get focusing on a prescribed point is easy to compute.
If $\frac{\sigma^2  L^3}{24 l_c} \geq r_0^2$, then $\alpha_L$ is given by (\ref{eq:phase2a})
and therefore one should know the statistics of the random medium  to get focusing on a prescribed point.
More exactly, if one wants to focus at the target point $(\bx_{\rm t},L)$, then one imposes the phase
\begin{equation}
\label{eq:phase2}
\bb_{\rm t} = \frac{1}{\alpha_L} \bx_{\rm t}.
\end{equation}
The shift that can be imposed is in fact limited by the SNR. Indeed, the peak intensity $
{\cal I}_{{\rm p},{\bf 0}}^\bb$ is damped when the shift becomes large:
\begin{eqnarray}
\nonumber
{\cal I}_{{\rm p},{\bf 0}}^\bb &=&
\lim_{\eps \to 0}  \Big| \EE \big[ \hat{u}_{\rm tr}^{\bb,\eps} (  \bx^\bb ; {\bf 0}) \big] \Big|^2\\
 &=& \frac{1}{ \big( 1 +  \frac{\sigma^2  L^3}{6 r_0^2 l_c} \big)^2} \exp \Big( - \frac{|\bb|^2}{2 r_0^2}  \big( 
 \frac{\frac{\sigma^2  L^3}{24 l_c}}{r_0^2 + \frac{\sigma^2  L^3}{24 l_c}} + \frac{\rho_0^2}{r_0^2-\rho_0^2} \big)\Big)  ,
\end{eqnarray}
while the mean intensity of the background fluctuations is:
\begin{eqnarray}
\nonumber
{\cal I}_{{\rm b},{\bf 0}}^\bb  &=&
 \lim_{\eps \to 0}  {\rm Var} \Big( \hat{u}_{\rm tr}^{\bb,\eps} ( \bx ; {\bf 0})  \Big)
 =
  \lim_{\eps \to 0}  \EE \Big[ \big| \hat{u}_{\rm tr}^{\bb,\eps} ( \bx ; {\bf 0}) \big|^2 \Big] 
  -  \Big| \EE \big[ \hat{u}_{\rm tr}^{\bb,\eps} ( \bx ; {\bf 0}) \big] \Big|^2  \\
  &=& \frac{1}{\big( 1 +  \frac{\sigma^2  L^3}{6 r_0^2l_c} \big)\big( 1 +  \frac{\sigma^2  L^3}{6 \rho_0^2 l_c}\big) } \exp \Big( - \frac{|\bb|^2}{2 r_0^2}  \frac{\rho_0^2}{r_0^2-\rho_0^2} \Big) ,
\end{eqnarray}
which is independent of $\bx$.
As a result the signal-to-noise ratio of the shifted peak at $\bx^\bb$ is
\begin{eqnarray}
\nonumber
{\rm SNR}^\bb_{\rm tr} &=&
\frac{{\cal I}_{{\rm p},{\bf 0}}^\bb }{{\cal I}_{{\rm b},{\bf 0}}^\bb } = \frac{ 1 +  \frac{\sigma^2  L^3}{6 \rho_0^2 l_c}}{  1 +  \frac{\sigma^2  L^3}{6 r_0^2 l_c}  }
 \exp \Big( - \frac{|\bb|^2}{2 r_0^2} \frac{\frac{\sigma^2  L^3}{24 l_c}}{r_0^2 + \frac{\sigma^2  L^3}{24 l_c}}\Big) \\
 &=&
{\rm SNR}_{\rm tr}
 \exp \Big( - \frac{|\bb|^2}{2 r_0^2} \frac{\frac{\sigma^2  L^3}{24 l_c}}{r_0^2 + \frac{\sigma^2  L^3}{24 l_c}}\Big) ,
\end{eqnarray}
where ${\rm SNR}_{\rm tr}$ is the signal-to-noise ratio (\ref{eq:snr1}) of the time-reversed refocused peak.
To observe refocusing, the ${\rm SNR}^\bb_{\rm tr}$ of the shifted peak should be larger than one,
and this means that one should limit the shift to $|\bb|\leq b_{\rm max}$ with
\begin{equation}
\label{def:bmax}
b_{\rm max}^2 =  2r_0^2 \Big(1+ \frac{24r_0^2l_c}{\sigma^2  L^3} \Big) \ln \frac{1 + \frac{\sigma^2  L^3}{6 \rho_0^2 l_c}}{1 + \frac{\sigma^2  L^3}{6 r_0^2 l_c}}   .
\end{equation}
This result with Eq.~(\ref{eq:phase2a}) shows that it is possible to focus in a region around the original source point whose radius $R_{\rm max} = \alpha_L b_{\rm max}$ is given by 
 \begin{eqnarray}
 \nonumber
  R_{\rm max}^2 &=& \frac{12l_c }{\sigma^2 k_0^2 L \big( 1 + \frac{\sigma^2  L^3}{24 r_0^2 l_c}\big)} \ln \frac{1 + \frac{\sigma^2  L^3}{6 \rho_0^2 l_c}}{1 + \frac{\sigma^2  L^3}{6 r_0^2 l_c}}  \\    
      &=&    \frac{    3 R^2_{\rm tr} }{ 1 + \frac{\sigma^2  L^3}{6 r_0^2 l_c}  }   
\ln {\rm SNR}_{\rm tr} 
 .
 \label{eq:Rmax1}
\end{eqnarray}
Thus, the focusing region increases as the logarithm of the ${\rm SNR}_{\rm tr}$ (\ref{eq:snr1}).

\section{Imaging through a Complex Medium}
\label{sec:image}%
\begin{figure}
\begin{center}
\begin{picture}(300,80)
\put(0,0){\includegraphics[width=6.1cm]{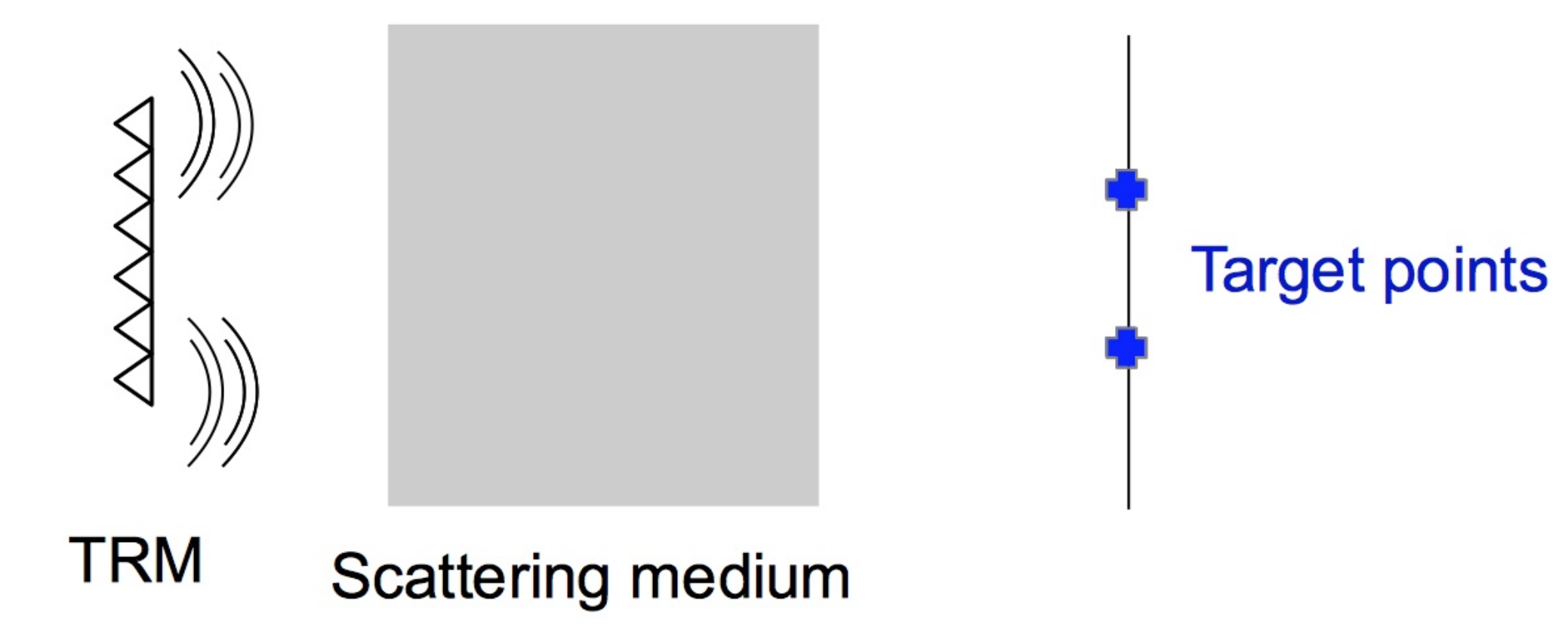}}
\put(170,-4){ \includegraphics[width=4.0cm]{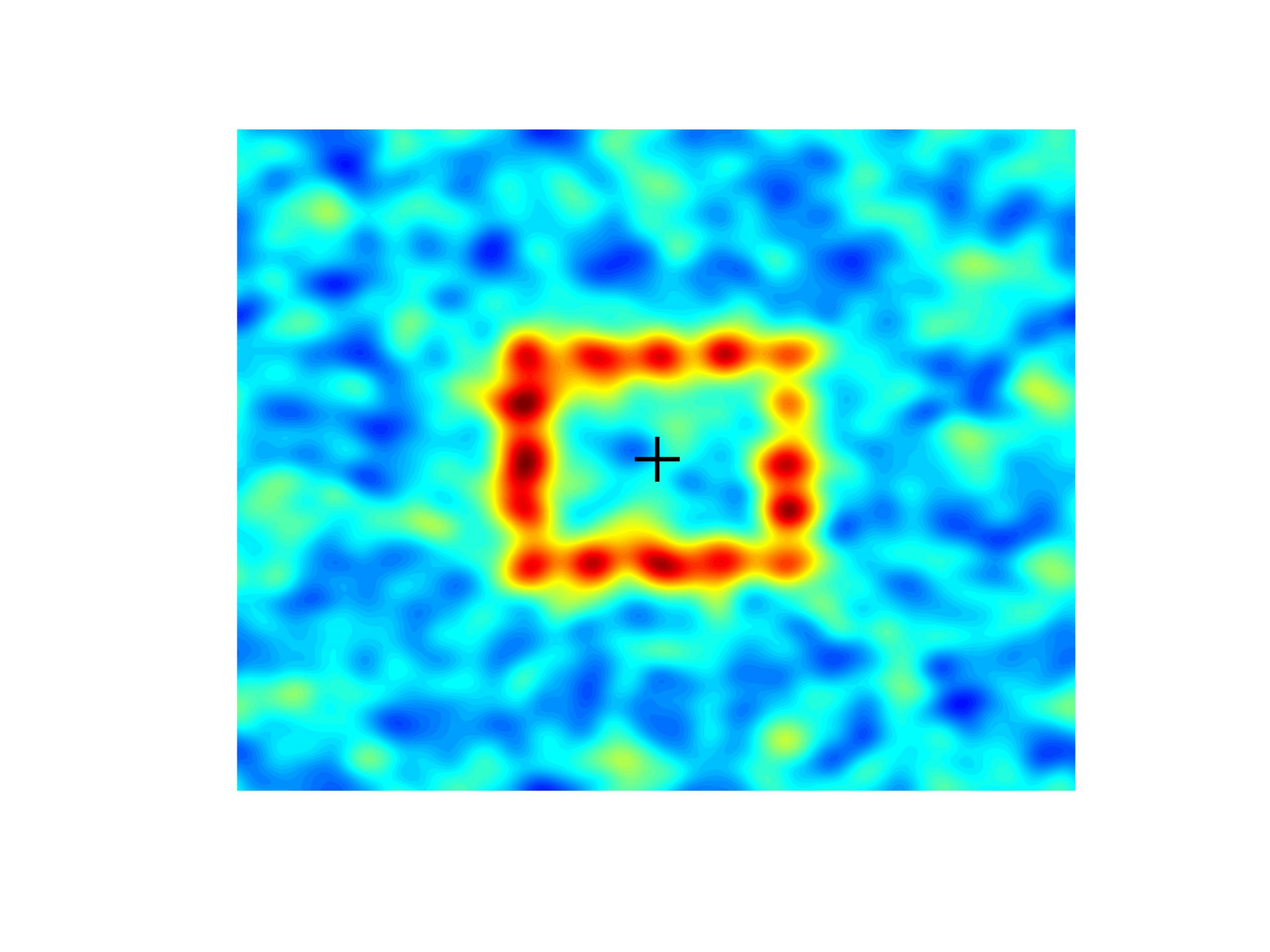}}
\end{picture}
\end{center}
\caption{Transmission of an image, here a square modeled as a set of sixteen target points.
The cross in the rigth image stands for the original source location.
\label{fig:0d} 
}
\end{figure}

The experimental observation that wave focusing can be achieved with the optimized 
phases of the SLM in the neighborhood of the target point and not only on the target point
is essential because it shows that an image can be transmitted once the phases
of the SLM have been optimized for a target point. This was achieved experimentally in 
\cite{katz12,popoff10} for instance.
In this section we quantify the resolution and stability properties of the transmitted image (see Figure \ref{fig:0d}).

In the time-reversal context, let $\psi(\bx)$ be a real-valued function that is the image to be transmitted. 
It could be a discrete image, i.e. $\psi$ is a sum of Dirac distributions, 
or a continuous one, i.e. $\psi$ is a bounded and compactly supported function.
Let us denote by $\hat{\psi}$ the Fourier transform of $\psi$:
$$
\hat{\psi} (\bxi ) =\int \psi(\bx) \exp ( - i \bx \cdot \bxi) d\bx .
$$
Let us emit the function $\overline{\hat{\psi} \big( {\bx_m}/{R_{\rm m}^2}  \big)}\overline{ \hat{u}_{\rm rec}(\bx_m;{\bf 0}) } $ 
from the  TRM and
 consider the refocused wave
\begin{equation}
\hat{u}_{\rm tr}^\psi (\bx; {\bf 0} ) =  \int 
\hat{u}_{\rm em} (\bx,\bx_m)
\exp \Big( - \frac{|\bx_m|^2}{R_{\rm m}^2} \Big) 
\overline{ \hat{\psi} \Big( \frac{\bx_m}{R_{\rm m}^2}  \Big) }
\overline{ \hat{u}_{\rm rec}(\bx_m;{\bf 0}) }  
  d\bx_m ,
\end{equation}
where $\hat{u}_{\rm rec}(\bx_m;{\bf 0})$ has been obtained with a point source at $({\bf 0},L)$.
By noting that the refocused wave can be expressed as
\begin{equation}
\hat{u}_{\rm tr}^\psi (\bx;{\bf 0} ) = \int 
\psi (\bb ) \hat{u}_{\rm tr}^\bb (\bx;{\bf 0}) d\bb ,
\end{equation}
with $\hat{u}_{\rm tr}^\bb $ defined by (\ref{eq:trfield2}),
the analysis of the image formation in the target plane $z=L$ follows from the results  obtained 
in the previous section.

In the scintillation regime, $C \to \eps C$, $r_0 \to r_0/\eps$, $\rho_0 \to \rho_0/\eps$, 
$L \to L/\eps$, $\psi(\bx) \to \eps^2 \psi(\eps \bx)$,
the refocused wave can be expressed as
\begin{equation}
\hat{u}_{\rm tr}^{\psi,\eps} (\bx;{\bf 0} ) = \int 
\psi(\bb) \hat{u}_{\rm tr}^{\bb,\eps} (\bx;{\bf 0}) d\bb .
\end{equation}
When scattering is strong $\sigma^2 k_0^2 l_c L \gg 1$ its expectation satisfies
\begin{eqnarray}
\nonumber
\EE \big[ \hat{u}_{\rm tr}^{\psi,\eps} 
 ( \bx  ;{\bf 0} )
 \big] &\stackrel{\eps \to 0}{\longrightarrow}&
\frac{1}{ 1 +  \frac{\sigma^2  L^3}{6 r_0^2 l_c} } 
 \int \frac{1}{\alpha_L^2} \psi
 \big(\frac{\bx'}{\alpha_L} \big)
 \exp \Big( - \frac{ |\bx-\bx'|^2}{2 R_{\rm tr}^2} \Big)\\
&& \times
 \exp \Big( - \frac{ |\bx'|^2}{4 r_0^2 \alpha_L^2} \big( 
 \frac{\frac{\sigma^2  L^3}{24 l_c}}{r_0^2 + \frac{\sigma^2  L^3}{24 l_c}} + \frac{\rho_0^2}{r_0^2-\rho_0^2} \big)\Big)d\bx' ,
 \label{eq:image1}
\end{eqnarray}
and its fluctuations are relatively smaller than its expectation provided the support of $\psi$
is within the disk with radius $b_{\rm max}$ defined by (\ref{def:bmax}). This gives an image of $\psi$
up to a dilatation by the factor $\alpha_L$ 
which is supported  within the disk with radius $R_{\rm max} = \alpha_L b_{\rm max}$.
The expression (\ref{eq:image1}) also shows that the dilated version $\psi(\cdot/\alpha_L)$ is imaged 
up to a radial attenuation and up to a convolution with a Gaussian kernel with radius $R_{\rm tr}$.
From (\ref{eq:Rmax1}) we have
$$
R_{\rm max}^2 =   \frac{ 3 R_{\rm tr}^2}{1+\frac{\sigma^2  L^3}{6 r_0^2 l_c}}  
\ln \frac{1 + \frac{\sigma^2  L^3}{6 \rho_0^2 l_c}}{1 + \frac{\sigma^2  L^3}{6 r_0^2 l_c}} .
$$
The transmitted image is acceptable if $R_{\rm tr} < R_{\rm max}$, otherwise the smoothing by the Gaussian kernel
with radius $R_{\rm tr}$ in (\ref{eq:image1}) blurs the image.
Therefore a favorable situation is when the radius of the TRM is large so that
$ {\sigma^2  L^3}/({6 r_0^2 l_c}) \lesssim 1$ and when it has a large number of elements 
so that $ {\sigma^2  L^3}/({6 \rho_0^2 l_c}) \gg 1$. 
In this situation $R_{\rm max} \gg R_{\rm tr}$ and the image is not significantly blurred, although
the propagation distance through the complex medium is much larger than the scattering mean free path.
This shows that imaging can be achieved through a strongly scattering medium by the SLM or
time-reversal technique.

\section{Conclusion}
Time reversal allows for wave refocusing through a complex medium.
This is true even when the elements of the TRM are larger than the 
correlation radius of the field that it records and that comes from a point source at a
target point on the other side of the medium. 
When refocusing on the target point, the profile of the mean focal spot depends on the diameter of the TRM but it 
does not depend on the number of elements of the TRM (provided it is large enough). However
 the signal-to-noise ratio strongly depends on number of elements of the TRM.
Moreover, we have shown that, when the field emitted by a point source in the target plane has been recorded by the TRM,
then it is also possible to focus a wave on a target point in the neighborhood of the original point source.
It is even possible to transmit an image.
The transmission of an image is possible provided the radius of the TRM is large enough and 
contains a large number of elements. All these results are quantified in this paper in the white-noise
paraxial regime, which is a regime relevant for laser beam propagation in scattering media and
in turbulent atmosphere in particular.

\section*{Acknowledgements}
We thank Dr. Arje Nachman for suggesting the above problem.  
This work is partly supported by AFOSR grant  \# FA9550-11-1-0176
and ANR project SURMITO.

\end{document}